\documentclass[12pt]{article}
\usepackage{amsfonts,amsmath,amsxtra}
\usepackage{amssymb,amscd}
\usepackage[mathscr]{eucal}
\def\hybrid{\topmargin 0pt      \oddsidemargin 0pt
        \headheight 0pt \headsep 0pt
        \textwidth 16.5cm
        \textheight 23cm
        \marginparwidth 0.0in
        \parskip 5pt plus 1pt   \jot = 1.5ex}
\catcode`\@=11
\def\marginnote#1{}
\newcount\hour
\newcount\minute
\newtoks\amorpm
\hour=\time\divide\hour by60 \minute=\time{\multiply\hour by60
\global\advance\minute by-\hour}
\edef\standardtime{{\ifnum\hour<12 \global\amorpm={am}%
        \else\global\amorpm={pm}\advance\hour by-12 \fi
        \ifnum\hour=0 \hour=12 \fi
      \number\hour:\ifnum\minute<10 0\fi\number\minute\the\amorpm}}
\edef\militarytime{\number\hour:\ifnum\minute<10 0\fi\number\minute}

\def\draftlabel#1{{\@bsphack\if@filesw {\let\thepage\relax
   \xdef\@gtempa{\write\@auxout{\string
      \newlabel{#1}{{\@currentlabel}{\thepage}}}}}\@gtempa
   \if@nobreak \ifvmode\nobreak\fi\fi\fi\@esphack}
        \gdef\@eqnlabel{#1}}
\def\@eqnlabel{}
\def\@vacuum{}
\def\draftmarginnote#1{\marginpar{\raggedright\scriptsize\tt#1}}

\def\draft{\oddsidemargin -0.1truein
        \def\@oddfoot{\sl preliminary draft \hfil
        \rm\thepage\hfil\sl\today\quad\militarytime}
        \let\@evenfoot\@oddfoot \overfullrule 3pt
        \let\label=\draftlabel
        \let\marginnote=\draftmarginnote
\def\@eqnnum{{\rm (\theequation)}
\rlap{\kern\marginparsep\tt\@eqnlabel}%
\global\let\@eqnlabel\@vacuum}  }

\newcommand{\RR}{{\mathbb{R}}}

\newfont{\Bbbb}{msbm7 scaled 1\@ptsize00}
\newcommand{\zs}{\raise-1pt\hbox{$\mbox{\Bbbb Z}$}}

scaled 1\@ptsize00

\font\sevenmsa=msam6 
\newfam\msafam
\textfont\msafam=\sevenmsa
\def\hexnumber@#1{\ifnum#1<10 \number#1\else
\ifnum#1=10 A\else\ifnum#1=11 B\else\ifnum#1=12 C\else \ifnum#1=13
D\else\ifnum#1=14 E\else\ifnum#1=15 F\fi\fi\fi\fi\fi\fi\fi}
\def\msa@{\hexnumber@\msafam}
\def\llcorner{\delimiter"4\msa@78\msa@78 }
\def\lrcorner{\delimiter"5\msa@79\msa@79 }
\mathchardef\blacktriangleright="3\msa@49
\mathchardef\blacktriangleleft="3\msa@4A \font\tenmsb=msbm10 scaled
1\@ptsize00
\newfam\msbfam
\textfont\msbfam=\tenmsb \scriptfont\msbfam=\tenmsb

\newdimen\linethick  \linethick=0.4pt
\newdimen\hboxitspace    \hboxitspace=5pt
\newdimen\vboxitspace    \vboxitspace=5pt

\def\fr#1{%
\be\new \vcenter{ \hrule height\linethick
           \hbox{\vrule width\linethick
                 \kern\hboxitspace
                 \vbox{\kern\vboxitspace
                       \hbox{$\begin{array}{c}\displaystyle#1
          \end{array}$}%
                       \kern\vboxitspace}%
                 \kern\hboxitspace
                 \vrule width\linethick}%
           \hrule height\linethick}%
\ee}

\newdimen\Squaresize \Squaresize=14pt
\newdimen\Thickness \Thickness=0.5pt

\def\Square#1{\hbox{\vrule width \Thickness
   \vbox to \Squaresize{\hrule height \Thickness\vss
      \hbox to \Squaresize{\hss#1\hss}
   \vss\hrule height\Thickness}
\unskip\vrule width \Thickness} \kern-\Thickness}

\def\Vsquare#1{\vbox{\Square{$#1$}}\kern-\Thickness}

\def\numberbysection{\@addtoreset{equation}{section}
        \def\theequation{\thesection.\arabic{equation}}}
\numberbysection

\renewcommand{\theequation}{\thesection.\arabic{equation}}
\def\titlepage{\@restonecolfalse\if@twocolumn\@restonecoltrue\onecolumn
     \else \newpage \fi \thispagestyle{empty}\c@page\z@
        \def\thefootnote{\fnsymbol{footnote}} }

\def\endtitlepage{\if@restonecol\twocolumn \else  \fi
        \def\thefootnote{\arabic{footnote}}
        \setcounter{footnote}{0}}  
\relax

\hybrid
\makeatletter
\newdimen\normalarrayskip            
\newdimen\minarrayskip               
\normalarrayskip\baselineskip \minarrayskip\jot
\newif\ifold             \oldtrue            \def\new{\oldfalse}
\def\arraymode{\ifold\relax\else\displaystyle\fi}
\def\eqnumphantom{\phantom{(\theequation)}} 
\def\@arrayskip{\ifold\baselineskip\z@\lineskip\z@
     \else
     \baselineskip\minarrayskip\lineskip1\baselineskip\fi}

\def\@arrayclassz{\ifcase \@lastchclass \@acolampacol \or
\@ampacol \or \or \or \@addamp \or
   \@acolampacol \or \@firstampfalse \@acol \fi
\edef\@preamble{\@preamble
  \ifcase \@chnum
     \hfil$\relax\arraymode\@sharp$\hfil
     \or $\relax\arraymode\@sharp$\hfil
     \or \hfil$\relax\arraymode\@sharp$\fi}}

\def\@array[#1]#2{\setbox\@arstrutbox=\hbox{\vrule
     height\arraystretch \ht\strutbox
     depth\arraystretch \dp\strutbox
width\z@}\@mkpream{#2}\edef\@preamble{\halign \noexpand\@halignto
\bgroup \tabskip\z@ \@arstrut \@preamble \tabskip\z@ \cr}%
\let\@startpbox\@@startpbox \let\@endpbox\@@endpbox
  \if #1t\vtop \else \if#1b\vbox \else \vcenter \fi\fi
  \bgroup \let\par\relax
  \let\@sharp##\let\protect\relax
  \@arrayskip\@preamble}
\def\eqnarray{\stepcounter{equation}%
              \let\@currentlabel=\theequation
              \global\@eqnswtrue
              \global\@eqcnt\z@
              \tabskip\@centering              
              \let\\=\@eqncr
              $$%
            \halign to \displaywidth  \bgroup
             \eqnumphantom \@eqnsel
      \hskip\@centering                               
    $\displaystyle  \tabskip\z@ {##}$%
    &\global\@eqcnt\@ne \hskip 2\arraycolsep
         $ \displaystyle  \arraymode{##}$\hfil
    &\global\@eqcnt\tw@ \hskip 2\arraycolsep
         $\displaystyle\tabskip\z@{##}$\hfil
         \tabskip\@centering
    &{##}\tabskip\z@\cr}
\makeatother

\def\IC{\mathbb{C}}

\def\IF{\mathbb{F}}

\def\IQ{\mathbb{Q}}
\def\IR{\mathbb{R}}
\def\IZ{\mathbb{Z}}
\def\CA {\mathcal{A}}

\def\CD {\mathcal{D}}

\def\CF {\mathcal{F}}

\def\CH {\mathcal{H}}

\def\CQ {\mathcal{Q}}
\def\CR {\mathcal{R}}
\def\CS {\mathcal{S}}
\def\CT {\mathcal{T}}
\def\CU {\mathcal{U}}

\def\a {{\alpha}}

\def\g {{\gamma}}

\def\l {{\lambda}}

\def\pr {\partial}

\def\zb {\bar{z}}

\def\Tr{{\rm Tr}}

\def\tt {{\vartheta}_{1}}
\def\frak{\mathfrak}
\def\Fg{{\frak g}}

\newtheorem{te}{Theorem}[section]
\newtheorem{de}{Definition}[section]
\newtheorem{prop}{Proposition}[section]           
\newtheorem{cor}{Corollary}[section]

\newtheorem{rem}{Remark}[section]
\newcommand{\beq}[1]{\begin{equation}\label{#1}}
\newcommand\eeq{\end{equation}}
\newcommand\bqa{\begin{eqnarray}}
\newcommand\eqa{\end{eqnarray}}
\def\be{\begin{eqnarray}\new\begin{array}{cc}}
\def\ee{\end{array}\end{eqnarray}}

\def\beq{\begin{equation}}
\def\eeq{\end{equation}}
\def\bse{\begin{subequations}}                
\def\ese{\end{subequations}}
\def\bp{\begin{pmatrix}}
\def\ep{\end{pmatrix}}

\def\i{\imath}

\newcommand{\R}{\mathbb{R}}
\newcommand{\C}{\mathbb{C}}

\newcommand{\F}{\mathbb{F}}

\newcommand{\Z}{\mathbb{Z}}

\def\stack#1#2{\raise0.7pt\hbox{$\mathrel{\mathop{#2}\limits^{#1}}$}}
\def\tr{\triangleright}
\def\tl{\triangleleft}
\def\sem{\mathsurround=0pt \raise1pt
\hbox{$\scriptscriptstyle>\!\!$}\:\!\!\tl}
\def\mes{\mathsurround=0pt \tr\!\:\!\raise0.8pt
\hbox{$\scriptscriptstyle\!\!<$}\,}
\def\]{\mathsurround=0pt ]\raise-2pt\hbox{$_\ast$}}

\def\la{\lambda}
\def\l{\lambda}

\def\<{\langle}
\def\>{\rangle}

\def\vk{\varkappa}
\def\CQ{{\cal Q}}
\def\frak{\mathfrak}

\def\CU{{\cal U}}

\def\CH{\mathcal{H}}

\def\we{\raise-1pt\hbox{$\,\stackrel{\wedge}{,}\,$}}
\def\tr{{\rm tr}\,}
\def\Tr{{\rm Tr}\,}
\def\pr {\partial}

\newcounter{pac}[section]

0
\begin{document}

\footnotesize
\normalsize

\thispagestyle{empty}

\begin{center}

\phantom.
\bigskip
{\hfill{\normalsize hep-th/yymmnnn}\\\
\hfill{\normalsize ITEP-TH-54/07}\\
\hfill{\normalsize HMI-07-09}\\
\hfill{\normalsize TCD-MATH-07-16}\\
[10mm]\Large\bf
On Baxter $\CQ$-operators \\ ${}$\\ And   Their Arithmetic
Implications\footnote{
Notes for ``Quantum Integrable Systems and Number Theory'' Seminar in
ITEP, Autumn 2007.}
 } \vspace{0.5cm}\\
\bigskip
\bigskip
{\large A. Gerasimov}
\\ \bigskip
{\it Institute for Theoretical and
Experimental Physics, 117259, Moscow,  Russia,} \\
{\it  School of Mathematics, Trinity
College, Dublin 2, Ireland, } \\
{\it  Hamilton
Mathematics Institute, TCD, Dublin 2, Ireland},\\
\bigskip
{\large D. Lebedev\footnote{E-mail: lebedev@itep.ru}}
\\ \bigskip
{\it Institute for Theoretical and Experimental Physics, 117259,
Moscow, Russia},\\
\bigskip
{\large S. Oblezin} \footnote{E-mail: Sergey.Oblezin@itep.ru}\\
\bigskip {\it Institute for Theoretical and Experimental Physics,
117259, Moscow, Russia}.\\
\end{center}

\vspace{0.5cm}

\begin{abstract}
\noindent

We consider Baxter $\CQ$-operators for
various versions of quantum affine Toda chain. The interpretation of
eigenvalues of the finite Toda chain Baxter operators as local Archimedean
$L$-functions proposed recently is generalized to the case of  affine Lie
algebras. We also introduce  a simple generalization of Baxter operators
and local $L$-functions compatible with this identification.
This gives a connection of the Toda chain Baxter $\CQ$-operators
with an Archimedean version of the Polya-Hilbert operator proposed by Berry-Kitting.
We also elucidate the Dorey-Tateo
spectral interpretation of eigenvalues of $\CQ$-operators.
Using explicit expressions for eigenfunctions of
affine/relativistic Toda chain
we obtain an  Archimedean analog of
Casselman-Shalika-Shintani formula for Whittaker function
in terms of characters.

\end{abstract}

\vspace{1cm}

\clearpage \newpage



\section{Introduction}

Theory of quantum integrable systems provides  a reach source
of new mathematical concepts. Probably the most famous example is the
notion of  quantum groups described in terms of $R$-matrices.
May be a less widely known is the notion of the $\mathcal{Q}$-operator
introduced by Baxter \cite{Ba}.  $\CQ$-operators appear  as
an important tool to solve quantum integrable systems.
In \cite{GLO} we introduce $\CQ$-operators for
$\mathfrak{gl}_{\ell+1}$-Toda chains (using previous results
of \cite{PG} for  periodic Toda chains) and
propose an interpretation of Baxter
integral $\CQ$-operators for $\mathfrak{gl}_{\ell+1}$-Toda chains
as particular elements of spherical Hecke algebra
$\CH(GL(\ell+1),\mathbb{R}), SO(\ell+1,\IR))$. Surprisingly
eigenvalues of the $\CQ$-operators acting on
class one principal series $\mathfrak{gl}_{\ell+1}$-Whittaker functions
are given by local Archimedean $L$-factors corresponding to the
principal series representations. This implies that the formalism of
$\CQ$-operators might shed light on a
hidden structure of the Archimedean ($\infty$-adic
arithmetic) geometry (see e.g. \cite{Man1} and references therein).

In this paper we discuss a set of generalizations of the results
presented in \cite{GLO}. First we introduce slightly generalized
$\CQ$-operators for finite Toda chins. This construction leads to a natural
connection between a representation of local Archimedean $L$-factors
(products of $\Gamma$-functions)  as  appropriately regularized
determinants of simple differential operators and the results of
Dorey-Tateo \cite{DT1,DT2} on a representation of eigenvalues of
$\CQ$-operators as spectral determinants of auxiliary differential
operators. The generalization of  \cite{DT1,DT2}  to Toda chains was
not known before.

Representations of $L$-factors
as determinants provide an intriguing relation with
the construction of the Weil group $W_{\IR}$ for
the field $\IR$ of real numbers. We demonstrate that the
operator arising in our interpretation of $\CQ$-operator's
eigenvalues for $\mathfrak{gl}_{\ell+1}$-Toda chain
is consistent with the structure of the Weil group $W_{\IR}$
and with the Archimedean analog of the Polya-Hilbert operator introduced
by Berry and Kitting \cite{BK}.

Following the interpretation of $\CQ$-operators for
$\mathfrak{gl}_{\ell+1}$-Toda chains given in \cite{GLO}
we consider  Toda chains for affine Lie algebras and
propose to identify the eigenvalues of the corresponding
$\CQ$-operators with the appropriate generalization of local
Archimedean $L$-factors for affine Lie algebras. We distinguish
affine Toda chains corresponding to  generic and critical
levels. The first one is known as Relativistic Toda chain \cite{Ru}
($q$-deformed Toda chain in the terminology of \cite{Et})
 and the second one is a
standard affine $\widehat{\mathfrak{gl}}_{\ell+1}$-Toda chain
proper. We construct explicitly $\CQ$-operators and
their eigenvalues for simplest non-trivial cases. Thus in the case
of $q$-deformed Toda chain an analog of $L$-factor is given by a
Double Sine function \cite{KLS}. To the best of our knowledge these
local Archimedean $L$-factors were not considered  before.

In this paper we  also consider a variant of a Toda chain
for affine Lie algebras at generic level with
eigenfunctions defined on a lattice. This type of Toda chain
naturally arises in the study of $K$-theory version of Gromow-Witten
invariants of flag varieties \cite{Gi,GiL}.  In the simplest case
we derive explicit formulas
for lattice analogs of  Whittaker functions and eigenvalues of
$\CQ$-operators.  Remarkably this provides an interpolation between
Archimedean and non-Archimedean  Whittaker
functions. There is a  well-known
identification of $p$-adic Whittaker functions for group $G(\IQ_p)$
and characters of finite-dimensional representations of the Langlands dual
complex Lie group ${}^LG_0(\IC)$  due to Casselman-Shalika-Shintani \cite{Sh,CS}.
We provide an expression for lattice Whittaker function as a trace
over a graded pieces of an infinite-dimensional representation of
a Heisenberg algebra. In the appropriate limit this leads to
an Archimedean variant of  Casselman-Shalika-Shintani formula.

The plan of this paper  is as follows. In Section 2
we propose  a generalization of the Baxter $\CQ$-operator for
$\mathfrak{gl}_{\ell+1}$-Toda chain and discuss a relations with
operators arising in \cite{DT1,DT2} and \cite{BK}. In Sections 3
and 4 we  consider various versions of Toda chains related to affine Lie
algebras,
corresponding $\CQ$-operators and the eigenfunctions (generalizations of Whittaker
functions).
Local archimedean $L$-factors of  affine Lie algebras
are proposed to coincide with the eigenvalues of the $\CQ$-operators.
In Section 5 an interpolation of Archimedean and non-Archimedean
 Whittaker functions and local $L$-factors is  constructed
using the a lattice version of the affine Toda eigenfunctions.
In Section 6 arithmetical interpretations of the constructions
presented  in the previous Sections are given. In particular the role
of the Weil  group $W_{\IR}$ is discussed and an Archimedean analog of
the Casselman-Shalika-Shintani  formula is proposed.

{\em Acknowledgments}: The research of AG was  partly supported by
SFI Research Frontier Programme and Marie Curie RTN Forces Universe
from EU. The research of SO was partially supported by  RF President
Grant MK-134.2007.1. The research was also partially supported by Grant
RFBR-08-01-00931-a.

\section{A Generalization of $\mathfrak{gl}_{\ell+1}$-Toda chain  $\CQ$-operator}

In this Section we propose a simple generalization of the Baxter
$\CQ$-operator for $\mathfrak{gl}_{\ell+1}$-Toda chain.
The  generalized $\CQ$-operators  reveals a connection between $\CQ$-operator
eigenvalues and spectral determinants of ordinary differential
operators. This is an analog of a spectral representation of Baxter
operator's eigenvalues in $c<1$ Conformal Field Theories due to
Dorey-Tateo \cite{DT1,DT2}.

We start by recalling the basic constructions in the theory
of quantum  $\mathfrak{gl}_{\ell+1}$-Toda chain.
Quantum $\mathfrak{gl}_{\ell+1}$-Toda chain  is a family
of mutually commuting quantum Hamiltonians generated by $(\ell+1)$
elementary Hamiltonians such that the first two are given by
\bqa \CH_1^{\mathfrak{gl}_{\ell+1}}&=&-\i\sum\limits_{i=1}^{\ell+1}
\frac{\partial}{\partial x_i},\\
\CH_2^{\mathfrak{gl}_{\ell+1}}&=&-\frac{1}{2}\sum\limits_{i=1}^{\ell+1}
\frac{\partial^2}{{\partial x_i}^2}+ \sum\limits_{i=1}^{\ell}
e^{x_{i}-x_{i+1}}.  \eqa
Common eigenfunctions of the Hamiltonians
can be identified with the $\mathfrak{gl}_{\ell+1}$-Whittaker functions \cite{Ko1,Ko2}
(see also \cite{Et})
 \bqa
\CH_1^{\mathfrak{gl}_{\ell+1}}\,\,\,\Psi^{\mathfrak{gl}_{\ell+1}}_{\g_1,\ldots
,\g_{\ell+1}} (x_1,\ldots,x_{{\ell}+1})&=&
\sum\limits_{i=1}^{\ell+1}\g_{i}\,\,\,
\Psi^{\mathfrak{gl}_{\ell+1}}_{\g_1,\ldots ,\g_{\ell+1}}
(x_1,\ldots,x_{\ell+1}),\\
\CH_2^{\mathfrak{gl}_{\ell+1}}\,\,\,\Psi^{\mathfrak{gl}_{\ell+1}}_{\g_1,\ldots
,\g_{\ell+1}} (x_1,\ldots,x_{{\ell}+1})&=&
\frac{1}{2}\sum\limits_{i=1}^{\ell+1}\g_{i}^2\,\,\,
\Psi^{\mathfrak{gl}_{\ell+1}}_{\g_1,\ldots ,\g_{\ell+1}}
(x_1,\ldots,x_{\ell+1}). \eqa

In \cite{GLO} the notion of the integral Baxter $\CQ$-operator for
$\mathfrak{gl}_{\ell+1}$-Toda chain  was introduced. Class one
principal series $\mathfrak{gl}_{\ell+1}$-Whittaker functions are
eigenfunctions of this operator with the eigenvalue given by the
corresponding local Archimedean $L$-factor. Let \be\label{renormW}
\Phi^{\mathfrak{gl}_{\ell+1}}_{\g_1,\ldots ,\g_{\ell+1}}
(x_1,\ldots,x_{\ell+1})=e^{\sum_{j=1}^{\ell+1}\rho_jx_j}\,
\Psi^{\mathfrak{gl}_{\ell+1}}_{\g_1,\ldots ,\g_{\ell+1}}
(x_1,\ldots,x_{\ell+1}), \ee where $\rho\in \RR^{\ell+1}$ , with
 $\rho_j=\frac{\ell}{2}+1-j,\,\,\,j=1,\ldots,\ell+1$.

\begin{te} Let $\CQ^{\mathfrak{gl}_{\ell+1}}(\l)$ be an  integral operator with the kernel
\be
\label{GLrator}
 \CQ^{\mathfrak{gl}_{\ell+1}}
(\underline{x},\underline{y}|\l)=2^{\ell+1}\exp\Big\{\sum_{j=1}^{\ell+1}(\i
\l +\rho_j) (x_j-y_j)-\\  -\pi
\sum_{k=1}^{\ell}\Big(e^{2(x_k-y_k)}+e^{2(y_k-x_{k+1})}\Big)- \pi
e^{2(x_{\ell+1}-y_{\ell+1})}\Big\}. \ee
Let $\Phi_{\g}^{\mathfrak{gl}_{\ell+1}}(x)$ be a modified
 $\mathfrak{gl}_{\ell+1}$-Whittaker function (\ref{renormW})
given by the integral representation
\bqa\label{giv}
\Phi_{\g_1,\ldots,\g_{\ell+1}}^{\mathfrak{gl}_{\ell+1}}
(x_1,\ldots,x_{\ell+1})=e^{\<\rho,x\>}\,
 \int_{\RR^{\frac{\ell(\ell+1)}{2}}}
\prod_{k=1}^{\ell}\prod_{i=1}^kdx_{k,i}\,\,
e^{\mathcal{F}^{\mathfrak{gl}_{\ell+1}}(x) }, \eqa where
\bqa\label{intrep}
\mathcal{F}^{{\mathfrak{gl}}_{\ell+1}}(x)=\imath\sum\limits_{k=1}^{\ell+1}
\g_k\Big(\sum\limits_{i=1}^{k}
x_{k,i}-\sum\limits_{i=1}^{k-1}x_{k-1,i}\Big)-
\pi\sum\limits_{k=1}^{\ell} \sum\limits_{i=1}^{k}
\Big(e^{2(x_{k+1,i}-x_{k,i})}+e^{2(x_{k,i}-x_{k+1,i+1})}\Big),\nonumber\eqa
and $x_i:=x_{\ell+1,i},\,\,\,i=1,\ldots,\ell+1$.
The following relation holds
 \be\label{BAxterredone}
\CQ^{\mathfrak{gl}_{\ell+1}}(\lambda)\, \cdot
\Phi_{\g_1,\ldots,\g_{\ell+1}}^{\mathfrak{gl}_{\ell+1}}
(x_1,\ldots,x_{\ell+1})=\prod_{j=1}^{\ell+1}\,
\pi^{-\frac{\imath\lambda-\imath\gamma_j}{2}}\,
\Gamma\Big(\frac{\imath \lambda-\imath \gamma_{j}}{2}\Big)
\,\,\Phi_{\g_1,\ldots,\g_{\ell+1}}^{\mathfrak{gl}_{\ell+1}}
(x_1,\ldots,x_{\ell+1}).\nonumber \ee
\end{te}
The eigenvalue of the $\CQ$-operator
\be\label{qEig}
 Q^{\mathfrak{gl}_{\ell+1}}(\lambda|\g_1,\ldots,\g_{\ell+1})\,
=\prod_{j=1}^{\ell+1}\,
\pi^{-\frac{\imath\lambda-\imath\gamma_j}{2}}\,
\Gamma\Big(\frac{\imath \lambda-\imath \gamma_{j}}{2}\Big), \ee is
identified with a local Archimedean $L$-factor corresponding to a
principal series representation induced from a Borel subgroup of
$GL(\ell+1,\IR)$ with a character parametrized by $(\g_1,\ldots
,\g_{\ell+1})$ (see \cite{GLO} for details). On the other hand the
eigenvalue (\ref{qEig}) is  an important ingredient of the Quantum
method of separation of variables \cite{Gu,Sk,KL1} for
$\mathfrak{sl}_{\ell+1}$-Toda chain. Quantum separation of variables
in a quantum mechanical system with $N$ degrees of freedom is a
unitary transformation relating Hamiltonian eigenfunctions of the
initial system with separated eigenfunctions of a new system given
by $N$ mutually non-interacting subsystems. Thus the separated
eigenfunctions of the transformed Hamiltonian are given by the
products of $N$ elementary eigenfunctions of auxiliary Hamiltonians.
Eigenfunctions  of $\mathfrak{sl}_{\ell+1}$-Toda chain can be
obtained from eigenfunctions
$\Phi_{\g_1,\ldots,\g_{\ell+1}}^{\mathfrak{gl}_{\ell+1}}
(x_1,\ldots,x_{\ell+1})$  by specializing spectral
variables $\sum_{j=1}^{\ell+1}\g_j=0$ and satisfy the invariance
condition
$$
\Phi_{\g_1,\ldots,\g_{\ell+1}}^{\mathfrak{sl}_{\ell+1}}
(x_1+a,\ldots,x_{\ell+1}+a)=\Phi_{\g_1,\ldots,\g_{\ell+1}}^{\mathfrak{sl}_{\ell+1}}
(x_1,\ldots,x_{\ell+1}).
$$
Quantum separation of variables in  $\mathfrak{sl}_{\ell+1}$-Toda chain
is given by the following integral
 \be
\Phi_{\g_1,\ldots,\g_{\ell+1}}^{\mathfrak{sl}_{\ell+1}}
(x_1,\ldots,x_{\ell+1})=\int_{\IR^{\ell+1}} \prod_{j=1}^{\ell}
dy_j\,
\tilde{\Phi}_{\g_1,\ldots,\g_{\ell+1}}^{\mathfrak{sl}_{\ell+1}}
(y_1,\ldots,y_{\ell})
\,\,R(x_1,\ldots,x_{\ell+1};y_1,\ldots,y_{\ell})\nonumber , \ee
where
$\tilde{\Phi}_{\g_1,\ldots,\g_{\ell+1}}^{\mathfrak{sl}_{\ell+1}}
(y_1,\ldots,y_{\ell})$ is a separated eigenfunction
$$
\tilde{\Phi}_{\g_1,\ldots,\g_{\ell+1}}^{\mathfrak{sl}_{\ell+1}}
(y_1,\ldots,y_{\ell})=\prod_{j=1}^{\ell}\phi(y_j;\g_1,\ldots
,\g_{\ell+1}),
$$
and
$$
 R(x_1,\ldots,x_{\ell+1};y_1,\ldots,y_{\ell})=\int
\prod_{j=1}^{\ell}\,d\tilde{\gamma}_j\,\,
 e^{-\i\sum_{j=1}^{\ell}\tilde{\gamma}_j\,y_j}
\,\,\Phi^{aux}_{\tilde{\gamma}}(x),
$$
$$
\Phi_{\tilde{\g}_1,\ldots,\tilde{\g}_{\ell}}^{aux}
(x_1,\ldots,x_{\ell+1})=e^{-\sum_{j=1}^{\ell}\tilde{\g}_j
x_{\ell+1}}\,\Phi_{\tilde{\g}_1,\ldots,\tilde{\g}_{\ell}}^{\mathfrak{gl}_{\ell}}
(x_1,\ldots,x_{\ell})
$$
is a kernel of the integral unitary transformation providing Quantum
separation of variables. Elementary eigenfunctions are given by
Fourier transform of the eigenvalue (\ref{qEig})
\be
\phi(y;\g_1,\ldots ,\g_{\ell+1})=\int_{\IR}\,\,d\l\,\, e^{\i \l y}\,
 Q^{\mathfrak{sl}_{\ell+1}}(\lambda|\g_1,\ldots,\g_{\ell+1}).
\ee The equations on the separated eigenfunctions
\be\label{finiteoperone}
(\prod_{j=1}^{\ell+1}(\i\gamma_j-\pr_y)-e^{-y})\phi(y;\g_1,\ldots
,\g_{\ell+1})=0, \ee is an instance of  the Baxter equation.  The
ordinary differential operator in (\ref{finiteoperone}) can be
considered as an  oper in the sense of \cite{BD} corresponding to
the  Whittaker function
$\Phi_{\g_1,\ldots,\g_{\ell+1}}^{\mathfrak{sl}_{\ell+1}}
(x_1,\ldots,x_{\ell+1})$ (on the relation between quantum separation
of variables and opers see \cite{Fr}). Note that this oper has
irregular singularities.

Rather surprisingly the integral operator (\ref{GLrator})
is non-trivial even for $\ell=0$ (the quantum separation of variables
is trivial in this case).
\begin{cor} The integral operator
$\CQ(\l)^{\mathfrak{gl}_{1}}$ with the kernel
\be
\CQ^{\mathfrak{gl}_{1}}
(x,y|\l)=2e^{\i\l (x-y)-\pi e^{2(x-y)}},
\ee
acts on the $\mathfrak{gl}_{1}$-Whittaker function
\bqa\label{givone}
\Phi_{\gamma}^{\mathfrak{gl}_1}
(x)= e^{\i \gamma x},
\eqa
as follows
\be\label{BAxterredonee}
\CQ^{\mathfrak{gl}_1}(\lambda)\,
\cdot\Phi^{\mathfrak{gl}_1}_{\gamma}
(x)=\pi^{-\frac{\imath\lambda-\imath\gamma}{2}}\,
\Gamma\Big(\frac{\imath \lambda-\imath \gamma}{2}\Big)
\,\,\Phi^{\mathfrak{gl}_{1}}_{\gamma}(x).
\ee
\end{cor}
The eigenvalue
\be
Q^{\mathfrak{gl}_1}(\l|\g)=\pi^{-\frac{\imath\lambda-\imath\gamma}{2}}\,
\Gamma\Big(\frac{\imath \lambda-\imath \gamma}{2}\Big),
\ee
is an elementary local Archimedean $L$-factor.

Now we introduce a generalized  Baxter $\CQ$-operator for the Toda chain
depending on an additional variable.
An  advantage of the generalized $\CQ$-operator is in
providing a simple relation with \cite{DT1,DT2}.

\begin{de} Define a generalized Baxter operator by the following
explicit formula for its integral kernel \bqa
\CQ^{\mathfrak{gl}_{\ell+1}}
(\underline{x},\underline{y}|\l,t)=2^{\ell+1}\exp\Big\{\sum_{j=1}^{\ell+1}(\i
\l +\rho_j) (x_j-y_j)-\\ \nonumber - t\pi
\sum_{k=1}^{\ell}\Big(e^{2(x_k-y_k)}+e^{2(y_k-x_{k+1})}\Big)- t\pi
e^{2(x_{\ell+1}-y_{\ell+1})}\Big\}.\eqa
\end{de}

\begin{prop}
  The following relations hold
\be \CQ^{\mathfrak{gl}_{\ell+1}}(\lambda,t)\,
\cdot\Phi^{\mathfrak{gl}_{\ell+1}}_{\underline{\gamma}}
(\underline{x})=Q^{\mathfrak{gl}_{\ell+1}}(\lambda,t)
\,\,\Phi^{\mathfrak{gl}_{\ell+1}}_{\underline{\gamma}}(\underline{x})=
\prod_{j=1}^{\ell+1}\, (\pi
t)^{-\frac{\imath\lambda-\imath\gamma_j}{2}}\,
\Gamma\Big(\frac{\imath \lambda-\imath \gamma_{j}}{2}\Big)
\,\,\Phi^{\mathfrak{gl}_{\ell+1}}_{\underline{\gamma}}(\underline{x}),\nonumber
\ee \be
 Q^{\mathfrak{gl}_{\ell+1}}(\lambda-2\i,t)=\Big(\prod_{j=1}^{\ell+1}\,
\frac{\imath\lambda-\imath\gamma_j} {2\pi
t}\,\Big)\,Q^{\mathfrak{gl}_{\ell+1}}(\lambda,t), \ee \be
 \Big\{\,t\frac{\pr}{\pr t}+\frac{1}{2}\sum_{j=1}^{\ell+1}\,
(\imath\lambda-\imath\gamma_j)\,\Big\}
\,Q^{\mathfrak{gl}_{\ell+1}}(\lambda,t)\,=\,0, \ee where
$\underline{\g}=(\g_1,\ldots , \g_{\ell+1})$,
$\underline{x}=(x_1,\ldots ,x_{\ell+1})$.
\end{prop}
It is useful to represent  $Q^{\mathfrak{gl}_{\ell+1}}(\l,t)$ as a determinant of a matrix
\be
Q(t,\l)=\det \widehat{Q}(t,\l),\qquad
\widehat{Q}(t,\l)=\Gamma(\frac{\i\l-\i\hat{\gamma}}{2}) \,(t\pi)^{
-\frac{\i\l-\i\hat{\gamma}}{2}},
\ee
 where $\hat{\gamma}={\rm diag}(\g_1,\ldots ,\g_{\ell+1})$ is a
 diagonal matrix.
Thus introduced  matrix $\widehat{Q}$ satisfies
difference-differential equations \be\label{Qteq} (t\frac{\pr}{\pr
t}+\frac{\i}{2}(\l-\widehat{\g}))\widehat{Q}(\lambda,t)=0, \qquad
\widehat{Q}(\l-2\i,t)=\frac{\i\l-\i\widehat{\g}}{2\pi
t}\widehat{Q}(\l,t). \ee

\begin{rem}
It  was argued in \cite{GLO} that eigenvalues of
$\CQ$-operators should be identified with local archimedean
$L$-factors. Given a generalization of $\CQ$-operator introduced  above
it is natural to consider a corresponding generalization of local
non-Archimedean $L$-factors
\be
L_p(s,n)=\int_{p^{n}\IZ_p}|x|^s\,d\mu_p(x)=\frac{p^{-ns}}{1-p^{-s}}.
\ee
Then we have
\be
L_p(s,n+1)=p^{-s}\,L_p(s,n),
\ee
and \be
\frac{L_p(s+1,n)}{L_p(s,n)}=p^{-n}\,\frac{L_p(s+1,0)}{L_p(s,0)}
,\qquad L_p(s,0)=L_p(s)=\frac{1}{1-p^{-s}}.
\ee
\end{rem}

The meaning of the matrix operator $\widehat{Q}$ is a fundamental
matrix of solutions of the  differential equation (\ref{Qteq}) for
$t\in [0,+\infty)$. The value at $t=1$ of the fundamental solution
of eigenvalue problem for a linear differential equation
$\CD(t,\pr_t)=t\frac{\pr}{\pr t}+\frac{\i}{2}(\l-\widehat{\g})$
 can be identified with the
spectral determinant of $\CD(t,\pr_t)$ (see \cite{DT1,DT2} for a
general discussion of this relation).  Thus we get an identification
of the spectral determinant of the operator $\CD(t,\pr_t)$
 with the Baxter $\CQ$-operator of
$\mathfrak{gl}_{\ell+1}$-Toda chain which is in agreement with the
results of \cite{DT1,DT2}. Note that this relation is equivalent to
the representation (\ref{Gdet}) defined below.

The spectral interpretation of the $\CQ$-operator's eigenvalues can
be supported by the following considerations (closely related with
that in \cite{C}). The following integral
representation for  $\log \Gamma(z)$ holds \be
\log\,\Gamma(\l-\g)=\int_0^{\infty} \frac{dt}{t} e^{-t}
\left((\l-\g-1)+\frac{e^{-(\l-\g-1)t}-1}{1-e^{-t}}\right)=\nonumber
\ee \be =\int_0^{\infty} \frac{dt}{t}
\,\left(\frac{e^{-(\l-\g)t}}{1-e^{-t}}+{\rm regulators}\right).
\label{Glogint}\ee

Consider the operator $K_t$ of scaling transformations
with weight $\g$ acting on the space of functions of one variable
\be
K_t\,f_\g(x)=e^{-t\g}f_{\g}(e^{-t}x). \ee
It can be represented as
an integral operator
\be K_t=e^{-t(\g+x\pr_x)},
\ee with the  kernel
\be
 K_t(x,y)=e^{-t\g}\delta(y-e^{-t}x).
\ee The trace of the operator is given by \be {\rm Tr}\,K_t=\int dx
e^{-t\g}\delta(x-e^{-t}x)=\frac{e^{-t\g}}{1-e^{-t}}. \ee Therefore
for appropriately regularized determinant (e.g using
$\zeta$-function regularization) \be\label{regdet} \log [{\rm
det}_{reg}(\g+x\pr_x)]^{-1}=-{\rm
  Tr}_{reg}\log(\g+x\pr_x)=\lim_{\epsilon\to 0}
{\rm
Tr}\Big(\int_{\epsilon}^{\infty}\frac{dt}{t}K_t-R(\epsilon)\Big),
\ee where $R(\epsilon)$ denotes singular terms such that   r.h.s. of
(\ref{regdet}) has a finite limit.
 Using the integral representation (\ref{Glogint}) we obtain
\be\label{Gdet}
 \Gamma(\l-\g)=[{\rm det}_{reg}(\l-\g+x\pr_x)]^{-1}.
\ee

Let us finally note that a very close determinant representation for
$\Gamma$-functions arises naturally in the representation theory of
a (double) of Yangian $Y(\mathfrak{gl}_1)$ \cite{KT}. Consider the
following representations of $Y(\mathfrak{gl}_1)$ (for the notations
see \cite{KT}) \be\label{Yrep}
\pi_1(h^-(z))=1-\frac{\g}{z}=\frac{z-\g}{z}, \qquad \qquad
\pi_2(h^+(z))=1-\frac{z}{\l}. \ee Then for the universal $R$-matrix
$\CR$ in the representations  (\ref{Yrep}) we have
\be\label{contfus} R_{\pi_1\otimes\pi_2}=(\pi_1(\l)\otimes
\pi_2(\g))\CR =\exp(\oint dz
\log\frac{\Gamma(z-\g)}{\Gamma(z)}\otimes \pr_z\log h^+(z))=\ee \be
=\frac{\Gamma(\l-\g)}{\Gamma(\l)}=[{\rm
det}_{reg}(1-\frac{\g}{\l+x\pr_x})]^{-1} .\nonumber \ee Note that we
can consider (\ref{contfus}) as a  continuous version of the fusion.
This defines the appropriate version of the representation of
$\CQ$-operator from the universal $R$-matrix
 for $\mathfrak{g}=\mathfrak{gl}_1$.

\section{Affine Toda chain:  generic level}

In this section we generalize some of the results of  \cite{GLO} to the
case of the Toda chains for affine Lie algebras. The construction is
based on the generalization \cite{Et}
 of the Kostant construction \cite{Ko1,Ko2} to affine Lie
algebras. We consider affine Toda chains corresponding to
representations of affine Lie algebras at a generic level and critical level separately.
The main objective is the identification of eigenvalues of analogs of  Baxter
$\CQ$-operators with the local Archimedean $L$-factors corresponding
to affine Lie groups. Therefore, taking into account the discussion in
the previous Section we mainly consider the first
non-trivial case of $\widehat{\mathfrak{gl}}_1$ and briefly
discuss generalizations to $\widehat{\mathfrak{gl}}_{\ell+1}$, $\ell>0$.

To introduce  notations, recall the standard definitions for affine
Lie algebra $\widehat{\mathfrak{gl}}_1$. Affine Lie algebra
$\widehat{\mathfrak{gl}}_1$ is generated by $K,\,d,\,\,u_n$,
$n\in\Z$ with  the relations: \bqa
[u_n,u_m]=-n\,K\delta_{n+m,0},\quad [K,u_n]=0,\quad [d,u_n]=nu_n.
\eqa We will consider the following representations of the universal
enveloping algebra $\CU(\widehat{\mathfrak{gl}}_1)$. A
representation $\pi_{\g,\mu,\vk}$ is generated by the highest vector
$|v_+\>$ and satisfies \be u_{n>0}|v_+\>=0,\quad u_0|v_+\>=\g
|v_+\>,\quad K|v_+\>=\vk|v_+\>. \ee This representation can be
realized in the space of polynomials of variables $x_n$, $n\in
\IZ_+$ and exponential functions of $x_0$ as
\bqa\pi_{\g,\,\vk}(u_{-k})=-k\,x_k,\hspace{1cm}
\pi_{\g,\,\vk}(u_0)=\frac{\pr}{\pr x_0}, \hspace{1cm}
\pi_{\g,\,\vk}(u_{k})=\vk\frac{\pr}{\pr\,x_k},\\ \nonumber
\pi_{\g,\,\vk}(K)=\vk,\hspace{1.5cm} \pi_{\g,\,\vk}(d)=\Big(\,\mu+
\frac{1}{2}\frac{\pr^2}{\pr x_0^2} -
\sum_{k=1}^\infty\,{k}x_k\frac{\pr}{\pr\,x_k}\,\Big).\eqa Another
type of representations is a generalization of the evaluation
representation to the case of the algebra with the derivative
operator $d$ \be\label{evrep} \pi_{{\rm ev}}(u_k)=\tau^{k},
\,\,\,\,\pi_{{\rm ev}}(d)=\tau\frac{\pr}{\pr\tau},\,\,\,\,\,
\pi_{{\rm ev}}(K)=0,\ee acting on the space of functions of $\tau$.
Left and right Whittaker vectors $\psi_{L}$, $\psi_R$ in the
representation $\pi_{\g,\vk}$ are defined by the conditions  \bqa
u_{n>0}\,\psi_R=0,\hspace{2cm}
u_{n<0}\,\psi_L=0,\hspace{2cm}u_0\psi_{R,\,L}=\i\g\,\psi_{R,\,L}.\eqa

Let us introduce a spherical subalgebra $\mathfrak{k}$ of
$\widehat{\mathfrak{gl}}_1$. Define an involution
$*\,:\,u_k\longmapsto -u_{-k}$.  The spherical
subalgebra $\mathfrak{k}\subset\widehat{\mathfrak{gl}}_1$ is defined
as a fixed point subalgebra with respect to the
involution\footnote{More generally for
  $\widehat{\mathfrak{gl}}_{\ell+1}$ we should consider
$g(\sigma)\to (g(2\pi -\sigma)^T)^{-1}$.}
 $*\,$. The spherical vector $\psi_K$ is an element of
the representation space of $\pi_{\g,\mu,\vk}$ defined by \bqa
(u_k-u_{-k})\psi_K=0,\qquad  n=1,2,\ldots. \eqa Explicitly
Whittaker and spherical vectors can be represented as follows
\be\psi_L=e^{\i\g x_0}\prod_{n=1}^\infty\delta(x_n),\qquad
\psi_R=e^{\i\g x_0},\qquad  \psi_K=e^{\i\g x_0}
\prod_{n=1}^\infty\,\sqrt{\frac{n}{2\pi\vk}}\,e^{-n\,x_n^2/2\vk}.\ee
There is a natural non-degenerate pairing between
$\CU(\widehat{\mathfrak{gl}}_1)$-modules spanned by $\psi_L$ and
$\psi_R$, and by $\psi_R$ and $\psi_K$ respectively: \be
\<\psi_L|\,\psi_R\>=
\int_{\R^{\Z_+}}\,\overline{\psi_L(x)}\,\psi_R(x)\,[dx]=1,\ee \be
\<\psi_K|\,\psi_R\>=\int_{\R^{\Z_+}}\,\,\overline{\psi_K(x)}\,\psi_R(x)\,[dx]=1,\qquad
[dx]=\prod_{n=1}^\infty\,dx_n.\nonumber\ee Let us parameterize
elements of the Cartan subgroup of $\widehat{GL}_1$ as \bqa
e^u=\exp\bigl\{\,t_0u_0+t_+K+t_-d\,\bigr\}.\eqa One defines
Whittaker function as: \be\label{Whittakeraf}
\Psi(t_0,\,t_-,\,t_+)\,=\,\<\psi_L|\,\pi_{\g,\vk,\,\mu}\bigl(
e^{t_0u_0\,+\,t_+K\,+\,t_-d})|\,\psi_R\>\,=\\
=\,\<\psi_K|\,\pi_{\g,\vk,\,\mu}\bigl(
e^{t_0u_0\,+\,t_+K\,+\,t_-d})|\,\psi_R\>\, =\exp\Big\{\i\g
t_0\,+\,t_+\vk\,+\, t_-\bigl(\mu-\frac{\g^2}{2}\bigr)\,\Big\}.\ee
The Whittaker function satisfies a set of equations given by the
projections of the functions of the generators of
$\widehat{\mathfrak{gl}}_1$ that commute with any element of
$\CU(\widehat{\mathfrak{gl}}_1)$. We consider only the part of
equations that govern the dependence on $t_0$. There is a general
construction of the central elements in terms of the spectral
invariants. In our case it is useful to consider the following
element of an appropriate generalization of
$\CU(\widehat{\mathfrak{gl}}_1)$ (we allow exponential functions of
generators, the standard  completion of
$\CU(\widehat{\mathfrak{gl}}_1)$  has a trivial center)
$$
\widehat{t}(\l)=[\det_{V}{}'\vk\,\,
\tau\pr_{\tau}]^{-1}\cdot\det_{V}(\vk\,\, \tau\pr_{\tau}+
 J(\tau)+\lambda),\qquad J(\tau)=\sum_{n\in \IZ}u_n\tau^n,
$$
where $\tau=\exp i \sigma$, $\sigma\in S^1=\IR/2\pi \IZ$,
$V=L^2(S^1)$ is the  space of square-integrable functions on the
circle $S^1$, and the determinant $\det'_{V}\vk\,\,
\tau\pr_{\tau}$ is taken over subspace orthogonal to the constant
functions on $S^1$.  Here we have used the generalized evaluation
representation (\ref{evrep}). The Whittaker function
(\ref{Whittakeraf}) satisfies the following equation
\bqa\label{Eigprop}
\widehat{t}(\l) \,\,\Psi^{\widehat{\mathfrak{gl}}_1}_\g(x)=
[\det{}'_{V}\vk\,\,
\tau\pr_{\tau}]^{-1}\,\det_{V}(\vk\,\, \tau\pr_{\tau}+\lambda-\g)\,\,\,
\Psi^{\widehat{\mathfrak{gl}}_1}_\g(x), \eqa where
$\Psi^{\widehat{\mathfrak{gl}}_1}_\g(x)=\Psi(t_0,t_+,t_-)$, $t_0=x$,
$t_{\pm}=0$. For the eigenvalue of $\widehat{t}(\l)$ we have the
following explicit expression (using the $\zeta$-regularization) \be
t(\l)=[\det{}'_{V}\vk\,\,
\tau\pr_{\tau}]^{-1}\cdot\det_V\bigl(\,\vk\tau\frac{\pr}{\pr\tau}+(\l-\g)\,\bigr)=\\
(\l-\g)\,\prod_{n\neq
0}\bigl(1-\i\frac{(\l-\g)}{n\vk}\,\bigr)=\frac{\vk}{\pi}\,{\rm sh}\,
 \pi\frac{\l-\g}{\vk}.
\ee Thus the eigenvalue problem  (\ref{Eigprop}) can be rewritten as
follows \be\label{gloneequation} (q^{\frac{\i}{2}
(\l+\i\pr_{x_0})}-q^{-\frac{\i}{2}(\l+\i \pr_{x_0})})\,
 \Psi^{\widehat{\mathfrak{gl}}_1}_\g(x)=
(q^{\frac{\i}{2} (\l-\g)}-q^{-\frac{\i}{2}(\l-\g)})
 \Psi^{\widehat{\mathfrak{gl}}_1}_\g(x),
\ee
where $q=\exp(2\pi \i/\vk)$.

Note that (\ref{gloneequation})
has infinite number of solutions because the difference operator in
r.h.s. commutes with any function of $x$ periodic with respect to
the shift $x\to x+2\pi\vk^{-1}$. Similarly to the case of the
wave-functions for finite Lie algebra Toda chains the choice of the
solutions is dictated by the analytic properties of the solutions.
As it was argued in \cite{KLS} the appropriate way
to impose the right analytic conditions is to use
additional ``dual'' equations
\be\label{gloneequationdual}
(\tilde{q}^{\frac{\i}{2} (\l+\i\pr_{x_0})}-\tilde{q}^{-\frac{\i}{2}(\l+\i \pr_{x_0})})\,
 \Psi^{\widehat{\mathfrak{gl}}_1}_\g(x)=
(\tilde{q}^{\frac{\i}{2} (\l-\g)}-\tilde{q}^{-\frac{\i}{2}(\l-\g)})
 \Psi^{\widehat{\mathfrak{gl}}_1}_\g(x),
\ee
where $\tilde{q}=\exp(2\pi \i\vk)$.

Using an analog of the Kostant reduction \cite{Ko1,Ko2} for affine
Lie algebra $\widehat{\mathfrak{g}}$ at the level $k$,  we obtain
Hamiltonians for relativistic ($q$-deformed) Toda chain
with $q=\exp(2\pi i/(k+h^{\vee}))$, $h^{\vee}$ being the dual Coxter
number. This is not surprising  taking into account a correspondence
between representations of $\widehat{\mathfrak{g}}$ of level $k$ and
$\CU_q(\mathfrak{g})$, $q=\exp(2\pi i/(k+h^{\vee}))$. Thus in the
case of $\widehat{\mathfrak{gl}}_{\ell+1}$-Toda chain the resulting
difference equations (analogs of  (\ref{gloneequation})) should be
eigenfunction equations for generating function of quantum
Hamiltonians of $q$-Toda chain corresponding to
$\mathfrak{gl}_{\ell+1}$. This theory was considered in the
framework of the quantum separation of variables in \cite{KLS}. In
the case of the $\mathfrak{sl}_{\ell+1}$-Toda chain we have \be
\Psi_{\g_1,\ldots,\g_{\ell+1}}^{\widehat{\mathfrak{sl}}_{\ell+1}}
(x_1,\ldots,x_{\ell+1})=\int_{\IR^{\ell}} \prod_{j=1}^{\ell}
dy_j\,
\tilde{\Psi}_{\g_1,\ldots,\g_{\ell+1}}^{\widehat{\mathfrak{sl}}_{\ell+1}}
(y_1,\ldots,y_{\ell})
\,\,R(x_1,\ldots,x_{\ell+1};y_1,\ldots,y_{\ell})\nonumber , \ee
where we imply that $\sum_{j=1}^{\ell+1}\g_j=0$ and
$\tilde{\Psi}_{\g_1,\ldots,\g_{\ell+1}}^{\widehat{\mathfrak{sl}}_{\ell+1}}
(y_1,\ldots,y_{\ell+1})$ is a separated eigenfunction
$$
\tilde{\Psi}_{\g_1,\ldots,\g_{\ell+1}}^{\widehat{\mathfrak{sl}}_{\ell+1}}
(y_1,\ldots,y_{\ell})=\prod_{j=1}^{\ell}\phi(y_j;\g_1,\ldots
,\g_{\ell+1}),
$$
and
$$
 R(x_1,\ldots,x_{\ell+1};y_1,\ldots,y_{\ell})=\int
\prod_{j=1}^{\ell}\,d\tilde{\gamma}_j \,\,\,
e^{-\i\sum_{j=1}^{\ell}\tilde{\gamma}_jy_j}
\widetilde{\Psi}^{aux}_{\tilde{\gamma}}(x),
$$
is a kernel of the
integral unitary transformation providing quantum separation of
variables for affine Toda chain. For the explicit form of the
integral kernel $\widetilde{\Psi}^{aux}_{\gamma}(x)$ see \cite{KLS}.
Elementary eigenfunctions are given by Fourier transform of the
eigenvalue (\ref{qEig}) \be \phi(y;\g_1,\ldots
,\g_{\ell+1})=\int_{\IR}\,\,d\l\,\, e^{\i \l y}\,
 Q^{\widehat{\mathfrak{sl}}_{\ell+1}}(\lambda|\g_1,\ldots,\g_{\ell+1}).
\ee The equations on the separated eigenfunctions
\be\label{finiteoper} (\prod_{j=1}^{\ell+1}(\imath
q^{\frac{\i}{2}(\g_j+\i\pr_y)}-\imath q^{-
\frac{\i}{2}(\g_j+\i\pr_y)})-e^{-y})\phi(y;\g_1,\ldots
,\g_{\ell+1})=0, \ee are equivalent to the  Baxter equation on
$Q^{\widehat{\mathfrak{sl}}_{\ell+1}}(\lambda|\g_1,\ldots,\g_{\ell+1})$.
The ordinary difference  operator in (\ref{finiteoper}) can be
considered as  a generalization  of the oper in the sense of
\cite{BD} corresponding to the  Whittaker function
$\Psi_{\g_1,\ldots,\g_{\ell+1}}^{\widehat{\mathfrak{sl}}_{\ell+1}}
(x_1,\ldots,x_{\ell+1})$. Note that this oper has irregular
singularities.

By an analogy with the finite Lie algebra case
we can expect that the Fourier transform of
an elementary wave-function $\phi(y;\g_1,\ldots
,\g_{\ell+1})$ entering the quantum separated problem provides
a reasonable candidate for the
eigenvalue of  the Baxter $\CQ$-operator for
$\widehat{\mathfrak{gl}}_1$. The Fourier transform
of the wave-function is given by
\bqa\label{Laf}
Q^{\widehat{\mathfrak{gl}}_1}(\i\l-\i \g)
=\sqrt{2\pi}\Big(\frac{2\pi}{\vk}\Big)^{\frac{1}{2}-\i\l+\i\g}
S_2^{-1}(\i\l-\i \g),
\eqa
where the function $S_2(z)$ is a specialization of the function
 defined in the Appendix B for $\omega_1=1$ and $\omega_2=\vk$.
The function $S_2^{-1}(z)$ satisfies the following functional
relations
\bqa
S_2^{-1}(z+1)=2\sin\frac{\pi z}{\vk}\cdot S_2^{-1}(z),\qquad
S_2^{-1}(z+\vk)=2\sin(\pi z)\cdot S_2^{-1}(z),
\eqa
that provide analogs of Baxter equation on the eigenvalue
$Q^{\widehat{\mathfrak{gl}}_1}(\i\l-\i \g)$. Note that when $\vk\to \infty$
the eigenvalues $Q^{\widehat{\mathfrak{gl}}_1}(\i\l-\i \g)$
reduces to the eigenvalues of the
Baxter operator for $\mathfrak{gl}_1$-Toda chain
\bqa
\lim_{\vk\rightarrow\infty}\,\,Q^{\widehat{\mathfrak{gl}}_1}(z)\,=\,
\Gamma(z).
\eqa
Introduce the function
\be
\CS(z)=e^{-\frac{\i\pi}{2}B_{22}(z)}\,S_2(z),\qquad
B_{2,2}(z)=\frac{1}{\vk}z^2-
\frac{1+\vk}{\vk}z+\frac{1+3\vk+\vk^2}{6\vk}.\nonumber
\ee

\begin{prop} Let   $\CQ^{\widehat{\mathfrak{gl}}_1}(\l)$
be an integral operator with the kernel
\be
\CQ^{\widehat{\mathfrak{gl}}_1}(x,y|\l)=
2\pi\sqrt{\vk}\int_{\IR-\i \epsilon} d\rho\,\,
e^{-\frac{\vk}{2\pi}\bigl(y-x +\ln\frac{2\pi}{\vk}-\rho-
\frac{\pi}{2}\frac{1+\vk}{\vk}\bigr)^2}\,e^{\rho(1+\vk)}\,
\CS^{-1}\bigl(-\i\frac{\rho\vk}{2\pi}\,\bigr),
\ee
where $\epsilon>0$ and $\vk>0$. Then the following holds
\be \CQ^{\widehat{\mathfrak{gl}}_1}(\l)\cdot
\Psi_{\g}^{\widehat{\mathfrak{gl}}_{1}}
(x)=\int \,dy \,\,\CQ^{\widehat{\mathfrak{gl}}_1}(x-y|\l)\,
\Psi_{\g}^{\widehat{\mathfrak{gl}}_{1}}(y)=
Q^{\widehat{\mathfrak{gl}}_1}(\i\l-\i\g)
\,\Psi_{\g}^{\widehat{\mathfrak{gl}}_1}(x).
\ee
\end{prop}
Thus $\CQ^{\widehat{\mathfrak{gl}}_1}(x,y|\l)$ is a kernel of the Baxter
$\CQ$-operator for $\widehat{\mathfrak{gl}}_1$. The proposed local
 Archimedean $L$-factor corresponding to a class of automorphic
representations of affine Lie groups is then given by
\bqa
L_{\infty}^{\widehat{\mathfrak{gl}}_{\ell+1}}(s|\g_1,\ldots
,\g_{\ell+1})=
\prod_{j=1}^{\ell+1} Q(\i s-\i \g_j),
\eqa
where $Q(s)$ is given by (\ref{Laf}).

\section{Affine Toda chain:  critical level}

In this Section we shortly comment on the
$\widehat{\mathfrak{gl}}_{\ell+1}$-Toda chain corresponding to
representations of affine Lie algebra
$\widehat{\mathfrak{gl}}_{\ell+1}$ at the critical level. In this
case there is a large center in the universal enveloping algebra
$\CU(\widehat{\mathfrak{gl}}_{\ell+1})|_{k=-h^{\vee}}$ \cite{FF}.
Linear and quadratic Hamiltonians of the corresponding Toda chain
are given by \bqa
\CH_1^{\widehat{\mathfrak{gl}}_{\ell+1}}(x_i)=-\i\sum_{i=1}^{\ell+1}
\frac{\partial}{\partial x_i},\eqa \bqa
\CH_2^{\widehat{\mathfrak{gl}}_{\ell+1}}(x_i)=-\frac{1}{2}\sum_{i=1}^{\ell+1}
\frac{\partial^2}{\partial x_i^2}+e^{x_{\ell+1}-x_1}+
\sum_{i=1}^{\ell}e^{x_{i}-x_{i+1}}. \eqa The corresponding integral
Baxter $\CQ$-operator was introduced in \cite{PG} \be
\CQ_0^{\widehat{\mathfrak{gl}}_{\ell+1}}
(x_1,\ldots,x_{\ell+1};\,y_1,\ldots,y_{\ell+1}|\l)=\exp\Big\{\,\sum_{j=1}^{\ell+1}
\imath \l (x_j-y_j)-
\sum_{i=1}^{\ell+1}(e^{x_i-y_i}+e^{y_{i+1}-x_i})\,\Big\},\nonumber\ee
where $y_{\ell+2}=y_1$.

In the case of $\widehat{\mathfrak{gl}}_{1}$-Toda theory one has \be
\CH_1^{\widehat{\mathfrak{gl}}_1}(x)=-\imath \frac{\pr}{\pr x},
\qquad \Psi^{\widehat{\mathfrak{gl}}_1}_{\g}(x)=e^{\i\g x}. \ee
Although the Hamiltonian $\CH_1^{\widehat{\mathfrak{gl}}_1}$ and
$\widehat{\mathfrak{gl}}_{1}$-Whittaker function are very simple the
corresponding  $\CQ$-operator is non-trivial \be
\CQ_0^{\widehat{\mathfrak{gl}}_{1}}
(x,y|\l)=\exp\Big\{\,\i\l(x-y)-e^{x-y}-e^{y-x}\,\Big\}, \ee
 \be
\CQ_0^{\widehat{\mathfrak{gl}}_{1}}
(\l)\cdot \Psi^{\widehat{\mathfrak{gl}}_1}_{\g}(x)
 =\int_{\IR}\, dy\,\,\CQ^{\widehat{\mathfrak{gl}}_{1}}
(x,y|\l) \Psi^{\widehat{\mathfrak{gl}}_1}_{\g}(y)=
Q_0^{\widehat{\mathfrak{gl}}_{1}} (\l,\g)\,
\Psi^{\widehat{\mathfrak{gl}}_1}_{\g}(x), \ee where the eigenvalue
$Q_0^{\widehat{\mathfrak{gl}}_{1}} (\l,\g)$ is given by the  value
of the normalized Macdonald function $2K_{\i\l-\i\g}(2e^{\tau})$ of
zero argument $\tau=0$ \be\label{Eigen}
Q^{\widehat{\mathfrak{gl}}_{1}}_0(\l,\g)
=2K_{\i\l-\i\g}(2)=\int\limits_0^{\infty}\frac{dv}{v}\,\,
v^{-\i(\l-\g)} e^{-(v+v^{-1})}, \ee where the following integral
formula for Macdonald function is used \be
K_{\nu}(2e^{\tau})=\frac{1}{2}\int_{0}^{\infty}v^{-\nu-1}e^{-e^{\tau}(v+{v}^{-1})}dv.\ee
The eigenvalue \eqref{Eigen} satisfies the Baxter equation
 \be ((\i\l-\i\g)+e^{\imath\pr_{\l}}-e^{-\imath\pr_{\l}})
Q_0^{\widehat{\mathfrak{gl}}_{1}} (\lambda,\g)=0. \ee Let us note
that this elementary Baxter operator is compatible with the quantum
separation of variables in affine Toda  chain for
$\widehat{\mathfrak{sl}}_{\ell+1}$ \cite{Gu,Sk,KL2} \be
\Psi^{\widehat{\mathfrak{sl}}_{\ell+1}}_{\g_1,\ldots \g_{\ell+1}}
(x_1,\ldots ,x_{\ell+1})=\int
d\underline{\tilde{\g}}\,\mu(\underline{\tilde{\g}})
\prod_{j=1}^{\ell}\,
Q_0^{\widehat{\mathfrak{sl}}_{\ell+1}}(\tilde{\g}_j,\g_1,\ldots
\g_{\ell+1})
\,\,\Psi^{\rm aux}_{\tilde{\g}_1,\ldots \tilde{\g}_{\ell}}(x_1,\ldots, x_{\ell+1}),\\
\mu(\underline{\tilde{\g}})=
\prod_{k>j}\,\bigl|\Gamma(\i\tilde{\g}_k-\i\tilde{\g}_j)\,\bigr|^{-2},
\nonumber \ee where $\underline{\tilde{\g}}=(\tilde{\g}_1,\ldots
,\tilde{\g}_{\ell})$ and the separated wave function satisfies the
Baxter equation \be
 (\prod_{j=1}^{\ell+1}(\i\lambda-\i\gamma_j)+
e^{\imath\pr_{\l}}+(-1)^{\ell+1}e^{-\imath\pr_{\l}})
Q_0^{\widehat{\mathfrak{sl}}_{\ell+1}}(\l,\g_1,\ldots,
\g_{\ell+1})=0. \ee In the coordinate representation we have
\be\label{operoned}
(\prod_{j=1}^{\ell+1}(\pr_y-\i\g_j)+e^{-y}+(-1)^{\ell+1}e^y)\psi(y,\g_1,\ldots,
\g_{\ell+1} )=0, \ee where \be
\psi(y,\g_1,\ldots\g_{\ell+1})=\int_{\IR}d\l\,\,e^{\i \l y}\,\,
Q_0^{\widehat{\mathfrak{sl}}_{\ell+1}}(\l,\g_1,\ldots \g_{\ell+1}).
\ee
This can be interpreted as $\widehat{\mathfrak{gl}}_{\ell+1}$-oper (in terminology of
\cite{BD}) corresponding to affine Whittaker function. Let us stress
that it has irregular singularities.

Note that (\ref{Eigen}) can be considered as an eigenfunction
$\Psi^{\mathfrak{sl}_2}_{\imath\l-\imath\g}(e^{\tau})$ of
$\mathfrak{sl}_2$-Toda chain at $\tau=0$. Thus it is natural to
define a generalized $\CQ$-operator for $\widehat{\mathfrak{gl}}_1$
such that its eigenvalues are given by eigenfunction
$\Psi^{\mathfrak{sl}_2}_{\imath\l-\imath\g}(t)\,=2\,K_{\imath(\lambda-\gamma)}(2e^{\tau}),
\,t=e^{\tau}$. The eigenvalues
$Q_0^{\widehat{\mathfrak{gl}}_1}(t,\l)$ of the generalized
$\CQ$-operator  satisfy the relations \be
Q_0^{\widehat{\mathfrak{gl}}_1}(t,\l+\i)-
Q_0^{\widehat{\mathfrak{gl}}_1}(t,\l-\i)+ \frac{\i\l-\i\g}{t}
Q_0^{\widehat{\mathfrak{gl}}_1}(t,\l)=0, \ee \be\label{auxoper}
(\pr_{\tau}^2+4e^{2\tau}+(\l-\g)^2)\,Q_0^{\widehat{\mathfrak{gl}}_1}(t,\l)=0
,\quad t=e^{\tau}. \ee Considered interpretation of the Baxter
$\CQ$-operator as a restriction of the wave function of the
$\mathfrak{sl}_2$-Toda chain at $\tau=0$ is consistent with
\cite{DT1,DT2}, where the same type of relation was established  for
the $\CQ$-operator describing integrable structure of $c< 1$
Conformal Field Theories. The arguments from \cite{DT1,DT2} imply
that the eigenvalue of $\CQ$-operator for
$\widehat{\mathfrak{gl}}_1$ is equal to  the spectral determinant of
the operator (\ref{auxoper}). The generalization of this spectral
interpretation for $\widehat{\mathfrak{gl}}_{\ell+1}$, $\ell\geq 1$
will be given elsewhere.

\section{On a lattice version of affine Whittaker functions}

In this subsection we discuss a different definition of the
Whittaker function for affine Lie group in the case  of generic
level. We will use a term $q$-deformed
Toda chain for the quantum integrable system
such that the affine Whittaker function is a corresponding common eigenvalue
(see \cite{Et} for a detailed description).  It is easy to see that a solution of
 $q$-deformed $\widehat{\mathfrak{gl}}_1$-Toda chain equation (\ref{Eigprop})
is unique when the eigenfunction $\Psi^{\widehat{\mathfrak{gl}}_1}_{\g}(x)$
is considered as a function on the lattice $2\pi \vk^{-1}\IZ\in \IR$
\be
\Psi^{\widehat{\mathfrak{gl}}_1}_{\g}(2\pi n\vk^{-1})=e^{\imath 2 \pi
  n \vk^{-1} \g}=q^{n\g},\qquad n\in \IZ.
\ee

An analog of the Baxter $\CQ$-operator in this
case  can be extracted from  the explicit quantum separation of
variables for $q$-version of $\mathfrak{gl}_{\ell+1}$-Toda chain.
To simplify considerations below we allow complex continuation
with respect to $\vk$.
The approach of  \cite{KLS} is easily adopted to the case of
 eigenfunctions on the lattice with the following result for
 eigenvalues of the Baxter $\CQ$-operator
\be\label{qeiglat}
Q^{\widehat{\mathfrak{gl}}_{\ell+1}}_{lattice}(t|t_1,\ldots,t_{\ell+1})=
\prod_{i=1}^{\ell+1}\Gamma_q(t t_i^{-1}),
\ee
where a $q$-analog of $\Gamma$-functions is given by
(see Appendix B)
\be
\Gamma_q(t)=\prod_{m=0}^{\infty}\frac{1}{(1-tq^m)}.
\ee
and we use the notations $t=q^s$, $t_i=q^{\gamma_i}$.
Eigenvalues (\ref{qeiglat}) are
reduced to eigenvalues of
$\CQ$-operators for $\mathfrak{gl}_{\ell+1}$ in the limit $t=q^s$, $t_i=q^{\g_i}$,
$q\to 1$ (after simple renormalization of $\CQ$-operator). There is an
analog of the Mellin-Barnes integral representations \cite{KL1,GKL}
for $q$-Whittaker functions in terms of  $\Gamma_q$-functions and Jackson integral.
In the following we derive another representation of the $q$-version
of $\mathfrak{gl}_2$-Whittaker function which has a close relation
with Casselman-Shalika-Shintani  formula for $p$-adic Whittaker function
and will play an important role in the next Section.

The eigenfunction problem for $q$-deformed
$\mathfrak{gl}_2$-Whittaker function on the lattice (common
eigenfunction of $q$-versions of $\mathfrak{gl}_2$-Toda chain
Hamiltonians)  can be written in the following form \cite{KLS}
\be\label{qTodagltwo}
t_1^{-1}t_2^{-1}\Psi^{(q)}_{t_1,t_2}(n-1)+(1-q^{n+1})\Psi^{(q)}_{t_1,t_2}(n+1)=
(t_1^{-1}+t_2^{-1})\Psi^{(q)}_{t_1,t_2}(n), \ee
 We relate this with
the functional equation for
 $q$-version of Baxter $\CQ$-operator using the following
expansion of the eigenvalue of $\CQ$-operator (see Appendix B)
\be\label{expand}
Q^{\widehat{\mathfrak{gl}}_{\ell+1}}_{lattice}(t|t_1,\ldots,t_{\ell+1})
=\sum_{n_1,\cdots,n_{\ell+1}=0}^{\infty} \frac{\prod_{j=1}^{\ell+1}
(tt_j^{-1})^{n_j}}
{\prod_{i=1}^{\ell+1}\prod_{j_i=1}^{n_i}(1-q^{j_i})}=\nonumber \ee
\be =\sum_{n=0}^{\infty}t^{n} \left(\sum_{n_1+\cdots +n_{\ell+1}=n}
\frac{\prod_{k=1}^{\ell+1}t_k^{-n_k}}
{\prod_{i=1}^{\ell+1}\prod_{j_i=1}^{n_i}(1-q^{j_i})}\right)=
\sum_{n=0}^{\infty}t^{n} \chi^{(q)}_n
(t_1,\cdots , t_{\ell+1}). \ee where $t_i=q^{\g_i}$ and we introduce
the following function\footnote{ In terms of non-commutative
variables $\{t_k\}$ such that $t_it_j=q\,t_jt_i$, $i<j$ the function
(\ref{qcharacter}) can be written in the following form
$$ \chi^{(q)}_n (t_1,\ldots ,
t_{\ell+1})= \frac{(t_1+ \cdots +t_{\ell+1})^n}{(n)_q!}.
$$} \be\label{qcharacter}
\chi^{(q)}_n(t_1,\ldots , t_{\ell+1})= \sum_{n_1+\cdots
+n_{\ell+1}=n}\frac{t_1^{-n_1}\cdots t_{\ell+1}^{-n_{\ell+1}}
}{(n_1)_q!\cdots (n_{\ell+1})_q!}, \ee
 where $(n)_q!=(1-q)...(1-q^n)$.  Consider the eigenvalue of
the $q$-version of the  $\CQ$-operator for $\mathfrak{gl}_2$ \be
Q^{\widehat{\mathfrak{gl}}_{2}}_{lattice}(t|t_1,t_2)= \Gamma_q(t
t_1^{-1})\,\Gamma_q(t t_2^{-1}). \ee It satisfies the following
functional  equation \be
Q^{\widehat{\mathfrak{gl}}_{2}}_{lattice}(tq|t_1,t_2)=
(1-tt^{-1}_1)(1-tt^{-1}_2)\,Q^{\widehat{\mathfrak{gl}}_{2}}_{lattice}(t|t_1,t_2),
\ee
which leads to the equation for coefficients
$\chi^{(q)}_n(t_1,t_2)$ in (\ref{expand}) \be
t_1^{-1}t_2^{-1}\chi^{(q)}_{n-1}(t_1,t_2)
+(1-q^{n+1})\chi^{(q)}_{n+1}(t_1,t_2)=
(t_1^{-1}+t_2^{-1})\chi^{(q)}_{n}(t_1,t_2). \ee
 Thus using (\ref{qTodagltwo}) we can identify
 \bqa\label{gltwoCS}
\Psi^{(q)}_{t_1,t_2}(n)=\chi^{(q)}_{n}(t_1,t_2).
 \eqa
Note that in the limit $q\to 0$ the r.h.s. reduces to the character of
the finite-dimensional representation $V_n=S^n\IC^2$ of $\mathfrak{gl}_2$. This is
a special case of more general formula expressing Whittaker function
in terms of characters. A well-known example of such type of
relation is the  Casselman-Shalika-Shintani formula for $p$-adic
Whittaker function for group $G(\IQ_p)$ in terms of the  character
of the Langlands dual group ${}^LG_0(\IC)$ \cite{Sh}, \cite{CS}. Indeed
(\ref{gltwoCS})  in the limit $q\to 0$, $t_1=p^{-y_1}$,
$t_2=p^{-y_2}$  turns into  the Casselman-Shalika-Shintani
formula which in the case of $GL(2,\IQ_p)$ is a consequence of a
simple identity \be
\frac{1}{(1-p^{-y_1})(1-p^{-y_2})}=\sum_{n=0}^{\infty} \chi_{n}(p^{-\hat{y}}).
\ee
Here  $\hat{y}={\rm diag}(y_1,y_2)$ is a diagonal matrix  and
$\chi_{n}(\hat{y})$ is a character of  the finite-dimensional
representation $S^n\IC^2$ given by  $n$-th symmetric power  of the
standard two-dimensional representation of $\mathfrak{gl}_2$.
In the limit $t_i=q^{\g_i}$, $t=q^s$, $q\to 1$ the $q$-deformed
Whittaker function $\Psi^{(q)}_{t_1,t_2}(n)$ reproduces  a
specialization $\Psi^{\mathfrak{gl}_2}_{\g_1,\g_2}(x,0)$
of the standard $\mathfrak{gl}_2$-Whittaker function
 and the expansion
(\ref{expand}) turns into the known relation between local
Archimedean $L$-factors and $\mathfrak{gl}_2$-Whittaker (see e.g.
\cite{Bu}) \bqa\label{intTWO}
\Gamma(s+\g_1)\Gamma(s+\g_2)=\int_{-\infty}^{+\infty}\,dx\,\, e^{\i
sx}\,\,\Psi^{\mathfrak{gl}_2}_{\g_1,\g_2}(x,0). \eqa
There is a generalization of the integral formulas (\ref{intTWO}) arising
naturally in applications of the Rankin-Selberg method to analytic
continuations of global $L$-functions
\bqa\label{orthGL}
 \int\limits_{\RR^{\ell}}
\prod_{j=1}^{\ell}d{x}_{\ell}\,\,
\overline{\Psi}^{\mathfrak{gl}_{\ell}}_{\g_1,\ldots \g_{\ell}}
(x_1,\ldots,x_{\ell}) \Psi^{\mathfrak{gl}_{\ell+1}}_{\la_1,\ldots
,\la_{\ell+1}} (x_1,\ldots,x_{\ell},0)
=\prod_{i=1}^{\ell+1}\prod_{k=1}^{\ell}
\Gamma(\i\la_{\ell+1,i}-\i\gamma_{\ell,k}).
 \eqa
 This  type of integral
representations was considered in \cite{St} (see also \cite{GLO},
Lemma 4.2). Its $p$-adic counterpart follows from  the Cauchy identity \bqa
\prod_{i=1}^{\ell+1}\prod_{j=1}^{\ell+1}\frac{1}{1-z_iw_j}=
\sum_{\Lambda} \,\,s_{\Lambda}(z)\,s_{\Lambda}(w), \eqa where the
sum goes over Young diagrams and $s_{\Lambda}(z)$ is  the Schur
polynomial corresponding to a diagram $\Lambda$.
Using the expansion of the appropriate product of $\Gamma_q$ and
iterative structure of (\ref{orthGL}) one can find an analog of
 (\ref{gltwoCS}) for a generic $q$-deformed
$\mathfrak{gl}_{\ell+1}$-Whittaker function. An interpretation  of the
identity (\ref{gltwoCS}) and its generalization to more general
Whittaker functions will be discussed in the next Section.

\section{Towards Arithmetic interpretations}

In this Section we put the results of the previous Sections in the
proper  perspective.  The natural framework is provided by
a geometry of the compactification  divisor at infinity  (corresponding
to the Archimedean valuation) of the arithmetic
scheme $Spec(\IZ)$. To describe the geometric structures in the vicinity
of the compactification divisor it is essential to use a particular
generalization of the standard Galois group
provided by the notion of the Weil group.
In the description of the Weil group of a number field
we mainly follow \cite{W}, \cite{De}, \cite{T} (see also \cite{ABV}).

The Galois group for $\IR$  is given by ${\rm Gal}(\IC/\IR)=\IZ_2$ and
a non-trivial generator corresponds to the complex conjugation.
 By standard arguments of abelian
class field theory the Galois group ${\rm Gal}(\IC/\IR)$
can  be identified with the group of unites $U_{\IR}=\IR^*/\IR_+=\pm 1$. However this
picture is oversimplified. There is a canonical action of the
multiplicative group $\IC^*$ on the complexified cohomology of
compact non-singular algebraic varieties over $\IC$ which in many
respects is analogous to the action of Galois group ${\rm
Gal}(\overline{\IF}_p/\IF_p)$ on the \'{e}tale cohomology of schemes
over $\IF_p$ (we do not consider a more general structure of weight
filtrations of Hodge structures on non-compact singular algebraic
varieties). To take into account this action one should consider
Weil groups of $\IR$ generalizing the  Galois
${\rm Gal}(\IC/\IR)$. The  Weil group $W_F$ of a number field $F$ should
satisfy the following defining properties \cite{T}. First, there
should exist a homomorphism with a dense image in the natural
 topology on ${\rm Gal}(\overline{F}/F)$
$$
\phi:\,\, W_F\rightarrow {\rm Gal}(\overline{F}/F).
$$
Second, for any Galois extension $E$ of $F$, there should be
an inclusion
$$
r:  E^*\rightarrow W^{ab}_E,
$$
such that the composition
$$
e\circ \phi:\,\,\,E^*\rightarrow {\rm Gal}^{ab}(\overline{F}/E)
$$
is a basic isomorphism of abelian class field theory. Here $\overline{F}$ is
an algebraic closure of $F$ and $\Gamma^{ab}=\Gamma/[\Gamma,\Gamma]$
is the  abelianization of $\Gamma$.
Thus for $F=\IC$, ${\rm Gal}(\overline{\IC}/\IC)$ is trivial,
$W_{\IC}=\IC^*$, $\phi$ is trivial and $r$ is the identity map.
For $F=\IR$, one has  ${\rm Gal}(\IC/\IR)=\IZ_2$.
The Weil  group $W_{\R}=\IC^*\cup j\IC^*$
is generated  by a copy of
 $\C^*$ and an element $j$, subjected  to the relations:
$$
j x j^{-1}=\overline{x}, \qquad  j^2=-1\in \IC^*,
$$
with the maps
$$
\phi:\,\,W_{\IR}\rightarrow {\rm Gal}(\IC/\IR),\qquad
\phi(x)=1,\,\,\,\,\,\phi(jx)=-1,\quad x\in \IC,
$$
$$
r:\,\,\IR^* \rightarrow W^{ab}_{\R},\qquad r(x)=x.
$$
Note that $W_{\IR}$ is non-abelian and
its abelianization is $W^{ab}_{\IR}=W_{\IR}/[W_{\IR},W_{\IR}]=\IR^*$.

In many problems the  Weil group $W_F$ plays the role of the Galois
group of $F$. Thus in the Langlands correspondence L-packages of
automorphic admissible supercaspidal
representations of $G(F)$ are classified by conjugacy
classes of continuous homomorphisms from the Weil group
into  Langlands dual group ${}^LG$. In the case of Archimedean fields we have the
following explicit description of the spaces of homomorphism
(Langlands parameters). Let ${}^LG_0(\IC)$ be a complex Lie group, Langlands-dual
to $G$ and ${}^L\mathfrak{g}_0 ={\rm Lie}({}^LG_0)$. Then the set of
 conjugacy classes of continuous homomorphisms $\phi$ from $\IC^*$
 into ${}^LG_0$ with a semi-simple image  can  be
identified with the set of pairs $(\lambda_1,\lambda_2)\in
{}^L\mathfrak{g}_0\times {}^L\mathfrak{g}_0$
of semi-simple elements, subjected to the
following  requirements:
$$
 [\lambda_1,\lambda_2]=0, \qquad \exp
(2\pi \imath \lambda_1)=\exp (2\pi \imath \lambda_2).
$$
The
homomorphism $\phi$ is then given by:
$$
\phi(e^t e^{i\theta})=
\exp (t(\lambda_1 +\lambda_2))\exp (i\theta
 (\lambda_1 -\lambda_2)), \qquad
 \phi(z)= z^{\lambda_1}\overline{z}^{\lambda_2},\qquad z=te^{\i
   \theta}.
$$
Langlands parameters $\phi$ for ${}^LG_0$, $F=\IR$   can  be identified
with the set of pairs $(y,\lambda)\in {}^LG_0\times {}^L\Fg_0$,
$\lambda$ being  semisimple, subjected to the  conditions:
$$
[\lambda,Ad_y(\lambda)]=0, \qquad
y^2=\exp(2\pi i \lambda).
$$
In particular, for ${}^LG_0(\IC)=GL_1(\IC)$  Langlands parameters for
$F=\IC$ are given by  the quasi-characters
$\chi_{s,m}(z)=|z|^sz^{m}$, $s\in \IC$,
$m\in \IZ$. Similarly for $\F=\IR$ the corresponding
quasi-characters are given by
$\chi_{s,\epsilon}(z)=|z|^sz^{(1-\epsilon)/2}$, $s\in \IC$,
$\epsilon=\pm 1$. Langlands parameters classify
irreducible representations entering the decomposition
of the space of functions on $GL_1$.
The correspondence between quasi-characters of the Weil group $W_{\IC}$
and irreducible representations of $GL_1(\IC)$ can be described as follows.
Let $D=x\frac{\pr}{\pr x}$ be a vector field on
a one-dimensional torus $\IC^*$. The following operator
$$
\Pi(z,\bar{z})=e^{2\pi i\,D\,\log z +2\pi i \overline{D} \log \bar{z} },
$$
acts on the space $\CF$ of functions on $GL_1(\IC)=\IC^*$. Decompose
$\CF$  as a sum of  eigenspaces  of
$\Pi(z,\zb)$.  This provides  a  decomposition on
one-dimensional irreducible representations
with respect to the action of $GL_1(\IC)$  by translations.
The eigenvalues of $\Pi(z,\zb)$ corresponding to
a one-dimensional irreducible representation
are given by elements of $Hom(\IC^*,GL_1)$
in accordance with Langlands correspondence.
This construction can be straightforwardly generalized to the case of
algebraic tori $G=GL_1\times \ldots \times GL_1$.
Let us remark that a similar construction arises in representation
theory of lattice vertex algebras in Conformal Field Theory. In
particular as a part of the vertex algebra structure a $A$-graded
Frobenius algebras
$\CA=\oplus_{\a \in A} \CA_{\Delta_{\a}}$, $\Delta_{\a}\in\IR$,
with a semisimple action of $W_{\IR}$ appears.
The  action of $W^{ab}_{\IR}$ on one-dimensional $\IR$-vector spaces is given by\be
W^{ab}_{\IR}: \,\CA\rightarrow \CA, \qquad
t:\,\,\CA_{\Delta_{\a}}\rightarrow
t^{\Delta_{\a}}\,\CA_{\Delta_{\a}},\qquad t\in \IR_+.
\ee

Now we are ready to discuss the constructions of the previous Sections
from the arithmetic point of view. Note that local non-ramified
non-Archimedean $L$-factors over  prime numbers $p$ are  given by
\be\label{Lnonarch}
L_p(s)=\frac{1}{1-g_p\,p^{-s}}=\exp(\,-\sum_{n=1}^{\infty}\,
\frac{1}{n}\,p^{-sn}\,\,g_p^n), \ee where $g_p=p^{y}$ is a value of
a character of ${\rm Gal}(\overline{\IF}_p/\IF_p)$ on  a  Frobenius
homomorphism  $Fr_p$. The close analogy between  Weil group
$W_{\IR}$ and Galois groups of non-archimedean fields leads to a
search of an  analog of (\ref{Lnonarch}) for
the Archimedean field $\IR$. It seems that  (\ref{Gdet}) and
(\ref{Glogint}) are  good candidates for the proper analogs of
representations of $L_p(s)$  in (\ref{Lnonarch}).
This  matches with the identification $W^{ab}_{\IR}=\IR^*$
discussed above. Indeed the following integral representation
(\ref{Glogint}) holds
\be\label{intGamma}
 \Gamma(\l-\g)\sim
\exp(\int_{\IR_+}\,d\mu(t)\,\, e^{-t(\l-\g)}),\qquad
d\mu(t)=(1-e^{-t})^{-1}\frac{dt}{t},
\ee
 where $e^{-t(\l-\g)}$ is a
character of the action of an element $t\in W^{ab}_{\IR}$ on a
one-dimensional module and the integration goes over elements of
$W^{ab}_{\IR}/\{{\pm 1}\}=\IR_+$. Note that the measure $d\mu(t)$ is
non-trivial. The integral representation \eqref{intGamma}
of the logarithm of $\Gamma$-function  arises naturally in the related context as a
contribution at Archimedean places in the ``Explicit formulas'' in
Number theory (see e.g.  a related discussion in \cite{C}).

The Weil group $W_{\IR}$ manifest in other representations of
local Archimedean $L$-factors.
The determinant representation \eqref{Gdet} for $\Gamma$-function
seems to be  very close analog of the standard representation
of non-Archimedean $L$-factor as an inverse determinant and is
compatible with the interpretation of $W^{ab}_{\IR}$ as
an abelianization of the generalized  Galois group of $\IR$. This
interpretation should be compared with other closely related
proposals. According to Polya and Hilbert it is natural to look for
an operator acting on a Hilbert space such that
its spectrum  is given by non-trivial zeros of
 Riemann $\zeta$-function. In \cite{BK} Barry and Kitting proposed an
Archimedean analog of this  hypothetical Polya-Hilbert operator
reproducing  analytic properties of the local Archimedean factor of
the global Riemann $\zeta$-function of $\overline{Spec(\IZ)}$. More
precisely  it was shown that a  quantum Hamiltonian operator
$\hat{H}=x\pr_x$ should describe a smooth part of the distribution
of non-trivial zeroes of $\zeta$-function captured by a local
Archimedean $L$-factor. Thus $\hat{H}$ can be interpreted as
a local Archimedean version of the Polya-Hilbert operator and
is identified with a generator of $W^{ab}_{\IR}$.
This seems to be in a complete correspondence with the
representation of the local Archimedean $L$-factor (expressed in
terms of $\Gamma$-function) as an inverse determinant\footnote{Let
us also remark that there is  another representation
of $\Gamma$-function as an inverse determinant of the operator with
a  discrete spectrum \cite{D1,D2} (see also \cite{Man2}).}
(\ref{Gdet}). Let us also remark that
the inverse determinant representation \eqref{Gdet}
for $\Gamma$-function is also compatible with Dorey-Tateo type
representation of the Baxter $\CQ$-operator and thus provides
a connection between  Dorey-Tateo differential operator $\CD$
 and  a generator of the Weil group $W_{\IR}^{ab}$.

To support the discussed interpretation of various quantum Toda chain
constructions in terms of Archimedean geometry we consider Archimedean
counterpart  of the Casselman-Shalika-Shintani formula \cite{Sh,CS}.
One of the manifestations of the Langlands duality over $\IQ_p$ is the
representation of the $p$-adic Whittaker
function for algebraic group $G(\IQ_p)$
 as a character of a finite-dimensional
representation of the  Langlands dual group ${}^LG_0$.
It was proposed by  Shintani \cite {Sh} for $GL_{\ell+1}$
and by  Casselman-Shalika \cite{CS} for all semi-simple Lie groups.
Let  $G(\IQ_p)=GL(\ell+1,\IQ_p)$, $\hat{g}_p(\g)$
 be a semisimple conjugacy class of ${\rm   diag}(p^{\g_1},\ldots, p^{\g_{\ell+1}})$
 in $GL(\ell+1,\IC)$ corresponding to the image of the
Frobenius morphism. Let  $V_{(n_1,n_2,\ldots ,n_{\ell+1})}$ be a
finite-dimensional irreducible representation of $GL(\ell+1,\IC)$
corresponding to a partition $(n_1\geq n_2 \geq \ldots \geq
n_{\ell+1})$. Then for the $p$-adic Whittaker function corresponding
to a principal series representation with the character
$\chi^{B}_{(p^{\g_1},\ldots, p^{\g_{\ell+1}})}(g)$
of a Borel subgroup $B\subset G$ the following
holds
\be\label{CasSh} W^{(p)}_{(\g_1,\ldots,\g_{\ell+1})}({\rm
diag}(p^{n_1},\ldots,
p^{n_{\ell+1}}))=\Tr_{V_{(n_1,\ldots,n_{\ell+1})}}\,\hat{g}_p(\g).
\ee
We propose to consider (\ref{gltwoCS}) and its generalizations as an
analog of (\ref{CasSh}). Indeed $q$-deformed Whittaker functions on
the lattice defined in the previous Section provide an interpolation
between Archimedean and non-Archimedean Whittaker functions. In the
non-Archimedean limit $q\to 0$, $t_1=p^{-y_1}$, $t_2=p^{-y_2}$
(\ref{gltwoCS}) reduces to a special case of (\ref{CasSh}). To
give  an interpretation of  Archimedean Whittaker functions  as
asymptotics of  characters it would be enough to represent
$q$-deformed Whittaker function on the lattice as a trace of some
operator compatible with (\ref{CasSh}) in the limit $q\to 0$.  We
have for  eigenvalues of the Baxter $\CQ$-operators in the
$q$-deformed  case
$$
Q^{\widehat{\mathfrak{gl}}_{\ell+1}}_{lattice}(t|t_1,\ldots,t_{\ell+1})
=\prod_{i=1}^{\ell+1}\Gamma_q(t t_i^{-1})=
\sum_{n_1,\cdots,n_{\ell+1}=0}^{\infty} \frac{\prod_{j=1}^{\ell+1}
(tt_j^{-1})^{n_j}}
{\prod_{i=1}^{\ell+1}\prod_{j_i=1}^{n_i}(1-q^{j_i})}.
$$
Consider an algebra $\CA$ with  a set of generators
$\{a^+_{n,i},a_{n,i}\}$, $n\in \IZ_+$, $i=1,\ldots, \ell+1$  and
relations
 $$
[a^+_{n,i},a_{m,j}]=\delta_{n,m}\delta_{i,j},\qquad
[a^+_{n,i},a^+_{m,j}]=[a_{n,i},a_{m,j}]=0,\qquad n,m\in \IZ_+,
\,\,\,i,j=1,\ldots ,\ell+1.
$$
Define the following operators
$$
L_0=\sum_{j=1}^{\ell+1}\sum_{n=1}^{\infty}na_{n,j}^+a_{n,j},\qquad
J^{(j)}_0=\sum_{n=1}^{\infty}\,a^+_{n,j}a_{n,j},\qquad
I_0=\sum_{j=1}^{\ell+1}\sum_{n=1}^{\infty}\,a^+_{n,j}a_{n,j}=
\sum_{j=1}^{\ell+1} J^{(j)}_0.
$$
Note that one can construct a representation of a version of
polynomial loop algebra using the
generators $\{a^+_{n,i},a_{n,i}\}$ with the finite-dimensional
subalgebra $\mathfrak{gl}_{\ell+1}$ generated by
$J_0^{ij}=\sum_{n=1}^{\infty}a^+_{n,i}a_{n,j}$.

A representation $W_{\CA}$ of $\CA$ can be realized in the space of polynomials of
$\{a_n\}$, $n=1,2,\ldots$  and we have
$$
Q^{\widehat{\mathfrak{gl}}_{\ell+1}}_{lattice}(t|t_1,\ldots,t_{\ell+1})
=\Tr_{W_{\CA}}\,q^{-L_0}\,t^{I_0}\,\,\prod_{j=1}^{\ell+1}\,t_j^{-J^{(j)}_0}.
$$
Therefore one can identify
$\chi^{(q)}_{n}(t_i)$ entering the  expansion
(\ref{expand}) with the characters of $I_0$-graded subspaces
$W^{[n]}_{\CA}=\{v\in W_{\CA}|I_0v=\,n\, v\}$
$$
\chi^{(q)}_{n}(t_i)
=\Tr_{W^{[n]}_{\CA}}\,q^{-L_0}\,\,\,\prod_{j=1}^{\ell+1}\,t_j^{-J^{(j)}_0}.
$$
Note that $I_0$ commutes with $J_0^{ij}$ and thus $W_{\CA}$ supports
the action of $\mathfrak{gl}_{\ell+1}$.  Taking into account the
identification of $\chi^{(q)}_{n}(t_1,t_2)$ with
$q$-Whittaker functions for $\mathfrak{g}=\mathfrak{gl}_2$
we have
$$
\Psi^{(q)}_{t_1,t_2}(n)=
\Tr_{W^{[n]}_{\CA}}\,q^{-L_0}\,\,\,t_1^{-J_0^1}\,t_2^{-J_0^2}.
$$
In the asymptotic limit of $q\to 1$ one recovers the local Archimedean
$L$-function. Note that the natural framework of
representation theory over archimedean fields should be
a kind of asymptotic geometry (also known as a tropical geometry).

\section{Conclusions and further directions}

In this paper we consider  Baxter $\CQ$-operators
for various generalization of Toda chains and relate them with local
Archimedean $L$-factors. The relation between $\CQ$-operators and
$L$-factors seems provide a consistent
approach to construction new candidates for local Archimedean
$L$-factors. Thus, for instance,  in this paper we propose  Archimedean $L$-factors
hypothetically corresponding to generic level automorphic representations of affine
Lie algebras. It would be interesting to define its  non-Archimedean
counterparts and construct the corresponding
 complete $\zeta$-function satisfying certain  functional equations
(see \cite{Kap1,Kap2} for a construction of
Hecke algebras over two dimensional fields $\IQ_p((t))$ defined by Parshin \cite{P}).
Hopefully  further advances in  understanding  the
connection between formalism of Baxter $\CQ$-operators and
$L$-functions corresponding to automorphic representations will be
useful for theories of quantum integrable systems and of automorphic
representations.

Let us comment  on a possible application of the discussed
constructions of affine Lie groups $\CQ$-operators
to  the study of hidden
underlying structures of  algebraic geometry over archimedean
fields.  The very existence of the $\CQ$-operator
for finite-dimensional Lie groups introduced in \cite{GLO}
can be traced back to the  construction of a spectral curve for
the corresponding integrable systems. In  contrast to the
case of affine Lie algebras the meaning of the spectral curve for
finite-dimensional Lie algebras is not so straightforward. The simplest  way to
uncover the spectral curve and Baxter $\CQ$-operator
for finite-dimensional Lie algebras is by a
degeneration of the corresponding constructions for affine Lie
algebras.  This seems to be an instance of a general phenomena
that constructions over real/complex numbers should be understood as
degenerations of the corresponding constructions for fields of higher
dimensions (see e.g. \cite{Ch} for  a proposal to consider
Yangian as a hidden symmetry of the representation theory
over  Archimedean fields).  More explicitly the relevance of the
higher-dimensional point of view manifests an attempt to find an
Archimedean version of Casselman-Shalika-Shintani \cite{Sh,CS}
formula. Using a version of the $q$-deformed $\mathfrak{gl}_{\ell+1}$-Toda
chain (related with generic affine Toda chain) we construct
generalized  Whittaker functions interpolating between
non-Archimedean and Archimedean Whittaker functions.
These  interpolation formulas  are closely related with the Macdonald
philosophy of symmetric polynomials interpolation \cite{Mac}.  One of the
consequences of the existence of the Whittaker function
 interpolation is an  Archimedean
version of the Casselman-Shalika-Shintani \cite{Sh,CS}
 relation between principal series
$p$-Whittaker function for $G(\IQ_p)$ and characters  of
finite-dimensional representations of the Langlands dual group ${}^LG_0$.
The Archimedean $\mathfrak{gl}_{\ell+1}$-Whittaker functions
are identified with asymptotic of characters of
fixed grade components of infinite-dimensional algebras.

The construction of the  Archimedean version of the
Casselman-Shalika-Shintani  formula leads to an interesting
opportunity explicit realize  $W_{\IR}$-modules and
to interpret the $W_{\IR}$-modules  using symmetries  of the
conjectural Archimedean topos $ \CT op_{\IR}$ responsible for the
properties of the algebraic varieties over archimedean fields
\cite{Be}. A suitable candidate for  appropriate cohomology theory
might  be found in the papers  \cite{GK,Gi,GiL} on mirror constructions
for flag manifolds of semisimple  Lie groups. Note that according to
\cite{GK,Gi,GiL} one expects that cohomological and $K$-theoretic
Gromow-Witten invariants of flag manifolds for a semisimple group
$G$ should be given in terms of eigenfunctions of
${}^L\mathfrak{g}_0$-Toda chains (i.e.
${}^L\mathfrak{g}_0$-Whittaker functions)  and its $q$-deformations
where ${}^L\mathfrak{g}_0={\rm Lie}({}^LG_0)$ and ${}^LG_0$ is a
Langlands dual to $G$. The recent advances in the construction of
Whittaker functions for classical Lie groups using a torification of
flag varieties \cite{GLO1} provides  a direct way to prove this
conjecture. Thus it would be quite reasonable to conjecture that the
Type $A$ topological strings on $G$-flag varieties should provide a
realizations of Weil group  modules for Archimedean version of
Casselman-Shalika-Shintani formula. The discussion of an  explicit
realization of this idea obviously deserves a separate publication
and we are going  to return to this fascinating subject in the
future.  However taking into account  the previous considerations it
looks plausible that the hidden underlying structure of algebraic
geometry over complex numbers is described by a two-dimensional
quantum field theory.

\section{Appendix A: Special functions}

The following integral representations for logarithm
of $\Gamma$-function hold (see e.g. \cite {WW})
$$
\log \Gamma(z)=\int_0^{\infty}
\frac{dt}{t} e^{-t} \left((z-1)+\frac{e^{(z-1)t}-1}{1-e^{-t}}\right),
$$
$$
\log \Gamma(z)=(z-\frac{1}{2})\log z-z+\frac{1}{2}\log 2\pi
+\int_0^{\infty}\left(
\frac{1}{2}-\frac{1}{t}+\frac{1}{e^{t}-1}\right)\frac{e^{-zt}}{t}dt.
$$
Let $\omega_1,\,\omega_2\in\R_+$. Define functions
$S_2(z)$ and $\CS(z|\omega)$ as
\bqa\label{DSineIntegral}
S_2(z)=e^{\frac{\i\pi}{2}B_{2,2}(z)}\,\CS(z|\omega)=
e^{\frac{\i\pi}{2}B_{2,2}(z)}\,\exp\int_{\R+\i0}\frac{e^{zt}}
{(e^{\omega_1t}-1)(e^{\omega_2t}-1)}\frac{dt}{t}=\nonumber\\
=\exp\frac{1}{2\pi\i}\int_C
\frac{\sinh\bigl(z-\frac{\omega_1+\omega_2}{2}\bigr)t}
{\sinh\frac{\omega_1}{2}t\cdot\sinh\frac{\omega_2}{2}t}\ln(-t)
\frac{dt}{t},\eqa
where
$$
B_{2,2}(z)=\frac{1}{\omega_1\omega_2}z^2-
\frac{\omega_1+\omega_2}{\omega_1\omega_2}z+
\frac{\omega_1^2+3\omega_1\omega_2+\omega_2^2}{6\omega_1\omega_2},
$$
 and the contour $C$ is given in \cite{KLS}.
The function $S_2(z)$ satisfies  the following functional
equations.
$$ S_2(z+\omega_1)\,=\,\frac{1}{2\sin\frac{\pi
    z}{\omega_2}}\,S_2(z),\qquad
S_2(z+\omega_2)=\frac{1}{2\sin\,\frac{\pi z}{\omega_1}}\,S_2(z).
$$
By a deformation of contours $C$ and $\R+\i0$ one can extend
the definition of $S_2(z)$ and $\CS(z|\omega)$ to the complex values of
$\omega_1,\,\omega_2$ with ${\rm Im}(\omega_1/\omega_2)>0$. Then the
following holds:
$$
\CS(z|\omega)\,=\,
\cfrac{\prod\limits_{m=0}^\infty\bigl(1-q^{2m}e(z/\omega_2)\bigr)}
{\prod\limits_{m=1}^\infty
\bigl(1-\widetilde{q}^{-2m}e(z/\omega_1)\bigr)},
$$
 where
$e(z)=e^{2\pi\i z}$, $q=e(\frac{\omega_1}{2\omega_2})$,
$\tilde{q}=e(\frac{\omega_2}{2\omega_1})$. Fourier transform of the
double sine function can be expressed in terms of double sine function
as follows
$$
\int_{C'}dt\,\CS(\omega_1+\omega_2-\i
t-a|\omega)
e^{2\pi\i\frac{zt}{\omega_1\omega_2}}\,=\,\sqrt{\omega_1\omega_2}\,
e^{-\frac{\pi\i}{2}B_{2,2}(0)}\, \CS^{-1}(\i z|\omega)\,
e^{-2\pi\frac{za}{\omega_1\omega_2}},
$$
where the contour $C'$ is given e.g. \cite{KLS}.

\section{Appendix B: q-Math}

In this Appendix we follow mainly \cite{CK}.
We use the following definition of $q$-numbers
$$
[n]_q=\frac{1-q^n}{1-q}.
$$
Traditionally $q$-Gamma function is defined as
$$
\tilde{\Gamma}_q(x)=
\frac{(q;q)_{\infty}}{(q^x;q)_{\infty}}\,(1-q)^{1-x},\,\,\,\,\,\,\,
 (a;q)_n=\prod_{i=0}^n(1-aq^i).
$$
We use the following modified $\Gamma_q$-function
$$
\Gamma_q(t)=\frac{1}{(t;q)_{\infty}},\qquad \Gamma_q(tq)=(1-t)\Gamma_q(t).
$$
In the limit $q\rightarrow 0$, $t\rightarrow p^{-x}$ we  have
$$
\Gamma^{(p)}(x)\equiv \lim_{q\rightarrow 0,t=p^{-x}}(q;q)_{\infty}\,
 \Gamma_q(t)=\frac{1}{1-p^{-x}}.
$$
In the  limit $t=q^x$, $q=\exp(2\pi i \hbar)$, $\hbar\rightarrow 0$
the $\Gamma_q$-function is reduced to the standard $\Gamma$-function
$$
 \lim_{\hbar \rightarrow 0} \,\,(q;q)_{\infty}\,\Gamma_q(q^x) \hbar^{x-1}=\Gamma(x).
$$
The $q$-analogs of the  exponential  function are  given by two functions
$$
e_q(t)=\sum_{n=0}^{\infty}
\frac{t^{n}}{(q;q)_n}=\frac{1}{(t,q)_{\infty}},\qquad
E_q(t)=\sum_{n=0}^{\infty} \frac{q^{\frac{n(n-1)}{2}}\,
  t^{n}}{(q;q)_n}=(-t,q)_{\infty},
$$
satisfying  the relation $e_q(t)E_q(-t)=1$. We have with our
definition of $\Gamma_q(t)$
$$
\Gamma_q(t)=\,e_q(t).
$$
The limits of both $e_q(t)$ and $E_q(t)$
for $t=q^x$, $q\rightarrow 1$ are equal to the standard
exponent. In the limit $q\to 0$, $t\to p^{-x}$
we have
$$
e_p(x)=\frac{1}{1-p^{-x}},\,\,\,\,\,E_p(x)=(1-p^{x}).
$$
Jackson integral is  define as
$$
\int_{0}^{\infty}f(x)d_qx=(1-q)\sum_{j=-\infty}^{\infty}q^jf(q^j).
$$
For $q$-Gamma function there is an analog of the standard  integral representation
\bqa
\Gamma_q(t)&=&\frac{1}{\prod_{m=0}^{\infty}(1-q^{m})}\
\int \chi_{t}(z) E_q(-qz) z^{-1}d_qz=\nonumber
\\&=&\,
\sum_{m=0}^{\infty} \,\frac{t^{m}}{\prod_{j=0}^{m}(1-q^{j})}=
\frac{1}{\prod_{m=0}^{\infty}(1-t q^{m})}, \nonumber \eqa where
$\chi_t(z)=t^{v(z)}$, $v(q^n)=n$. In the limit $t=p^{-x}$, $q\to 0$
the integral expression is reduced to following
$$
\Gamma^{(p)}(x)=\sum_{j=0}^{\infty}p^{-xj}=\frac{1}{1-p^{-x}}.
$$

\end{document}